\documentclass[11pt]{article}

\usepackage{color,graphicx,times}
\usepackage{hyperref}

\newtheorem{lemma}{Lemma}

\newtheorem{theorem}{Theorem}
\newtheorem{corollary}{Corollary}

\newcommand{\qed}{\hfill {\bf q.e.d} }
\newcommand{\Bbb}[1]{%
{\bf #1 }}
\newcommand{\scr}[1]{%
\mathcal{#1}}
\newcommand{\di}[1]{%
\mathcal{#1}}
\newcommand{\vek}[1]{%
{\bf #1}}
\newcommand{\mat}[1]{%
\hat{#1}}
\newcommand{\sr}[1]{%
\mbox{{\scriptsize #1}}}

\newcommand{\proof}{ {\bf Proof: }}
\newcommand{\gen}{ {g}}
\newcommand{\osa}{ {g_{o}}}
\newcommand{\regg}{ {g_{r}}}

\newcommand{\op}{octahedral buddies }
\newcommand{\opc}{octahedral buddies}
\newcommand{\Op}{Octahedral buddies }
\newcommand{\OP}{Octahedral Buddies }

\newcommand{\K}{\mathrm{K}}
\newcommand{\wacky}{22.5\K}
\newcommand{\wack}{90\K}
\newcommand{\trivial}{K_4}
\newcommand{\negg}{Neg}

\newcommand{\vol}{ \mathcal{V}}

\begin{document}

\title{23040 symmetries of hyperbolic tetrahedra}

\author{Peter Doyle \and  Gregory Leibon}

\maketitle

\begin{abstract}

We give a rigorous geometric proof of the Murakami-Yano formula
for the volume of a hyperbolic tetrahedron.
In doing so, we are led to consider \emph{generalized hyperbolic tetrahedra},
which are allowed to be non-convex, and have vertices `beyond infinity';
and we uncover a group, which we call $\wacky$, of
$23040 = 30 \cdot 12 \cdot 2^6$
scissors-class-preserving symmetries of the space of 
(suitably decorated)
generalized
hyperbolic tetrahedra.
The group $\wacky$ contains the Regge symmetries as a subgroup of order
$144 = 12 \cdot 12$.
From a generic tetrahedron,
$\wacky$ produces 30 distinct generalized tetrahedra in the same scissors class,
including the 12 honest-to-goodness tetrahedra produced by the Regge subgroup.
The action of $\wacky$ leads us to the Murakami-Yano formula, and to 9 others,
which are similar but less symmetrical.
From here, we can derive yet other volume formulas with pleasant algebraic
and analytical properties.
The key to understanding all this is a natural
relationship between a hyperbolic tetrahedron and a pair of ideal
hyperbolic octahedra.

\end{abstract}

\begin{section}{Introduction \label{int}}

The computation and understanding of hyperbolic volume is an old and 
difficult problem.  Of fundamental interest and importance has been the
exploration of the volume of the hyperbolic tetrahedron. In particular,
historically there has been great interest in finding formulas for the
tetrahedron's volume which have optimal algebraic simplicity, together
with a concrete geometric interpretation (see Milnor \cite{Mi1},
Kellerhals  \cite{Ke}, and the references therein).

For the special case of an ideal hyperbolic tetrahedron,
Milnor \cite{Mi1} presents and derives, in a
straight-forward geometric way, what we presume is the optimally
elegant volume formula.
(See section \ref{itet} below.)

For a general hyperbolic tetrahedron,
Murakami and Yano
\cite{Mu}
recently found what we consider the most elegant known volume formula.
(See section \ref{MY} below.)
This formula arose from attempts to resolve
Kashaev's conjecture that the colored Jones polynomials
of a hyperbolic knot determines the hyperbolic volume of 
the knot's complement.
It was discovered utilizing properties
of the quantum $6j$-symbols,
and justified by means of a  known  formula for hyperbolic
volume due to Cho and Kim (see \cite{Ch}).
Murakami-Yano's derivation of their  formula was formal and was lacking a
concrete geometric interpretation.   
One goal of this paper is to provide  a rigorous 
geometric  interpretation, and explore
several new views of tetrahedral volume.

\begin{subsection}{ Volume Formulas \label{volfor}}

Most of the formulas in this paper will be described as   the volume of a specified scissors congruence class, hence we will first recall this concept.

\begin{subsubsection}{Scissors Congruence  \label{siscon}}

In this section we review the concept of scissors congruence and fix the notation that we will use  in order  to describe  a scissors class. 
To articulate the notion of scissors congruence needed in this paper,  we first form  the free Abelian group,  $\scr{F}$,   generated by the symbols $P$,  one for each  (unoriented) geodesic polyhedron $P \subset H^3$. 
Let 
\[\scr{R} =  \left<P+Q-R,\sigma(S)-S \right>, \] 
where the geodesic polyhedra $P$ and $Q$ have an intersection with empty interior and a union with interior equal to the interior of  the geodesic polyhedron $R$,  
and  where $\sigma$ is an orientation preserving isometry
being applied to  a  geodesic polyhedron $S$.   
The {\it scissors congruence group}, $\scr{P}(H^3)$,  is isomorphic to the quotient group  $ \scr{F}/\scr{R}$.   
 Notice that the volume, extended to $\scr{F}$ by linearity, provides a well defined homomorphism  
 \[ \vol:  \scr{P}(H^3) \rightarrow {\bf R}.  \]
   We call any pair of elements in $\scr{F}$ that agree in $\scr{P}(H^3)$ 
 {\it scissors congruent}.  Given a geodesic polyhedron  $P$ we will let $[P]^s$ denote its scissors class.
 
 There are a pair of observations about scissors classes that will be useful in what follows.
  First,  Dupont and Sah proved that we may divide by $2$ in  $\scr{P}(H^3)$ (see \cite{Sa1}).
In other words,  
 \begin{eqnarray} \label{half}
   &    
 \begin{array}{ccc}
& \mbox{if} & 2[P]=2[Q], \\
& \mbox{then} & [P]=[Q].  
\end{array}   &
 \end{eqnarray}  

The second concerns a  geodesic polyhedron's mirror image. 
 Given a  geodesic polyhedron $P$ let   $P^*$  denote its mirror image.  
 We have that 
 \begin{eqnarray}\label{mirror}
 [P]^s & = & [P^*]^s,
 \end{eqnarray} 
   as  was noted   in a letter to Gauss from  Gerling  in 1844 (see Neumann \cite{Ne}).

Notice that the generators of $\scr{F}$ have no orientation associated to them, hence the scissors class $[P]^s$ will ignore any orientation data associated to $P$. 
If we are given an oriented convex polyhedron $P$, then we will let
$[P]=[P]^s$ if $P$ is positively oriented and let   $[P]=-[P]^s$ if  $P$ is negatively oriented.
We extend this to lists of oriented convex polyhedra by letting  by letting 
\[  \left[(P_1,\ldots , P_M)\right] = \sum_{i=1}^M  [P_i]. \]

  \end{subsubsection}

\begin{subsubsection}{ The  Ideal Tetrahedron \label{itet}}

\begin{figure}
\center{\includegraphics{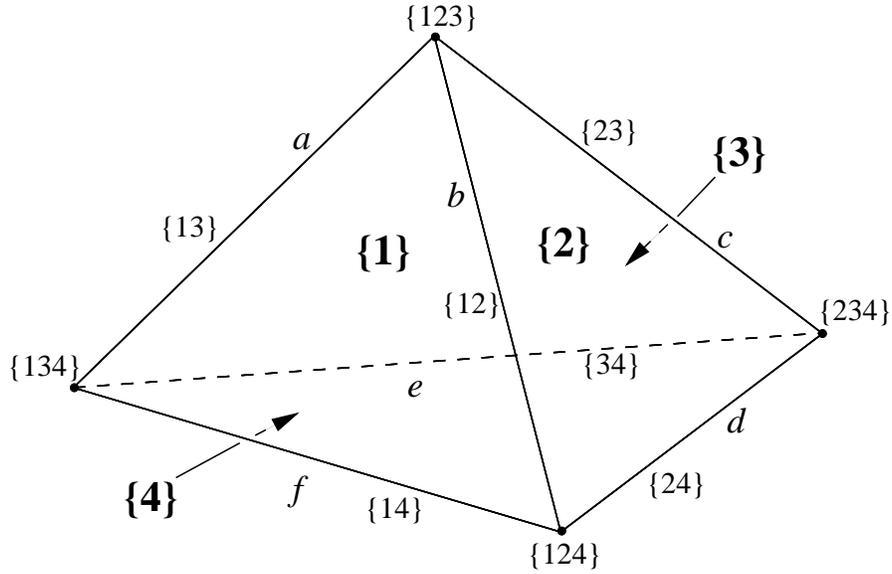}}
\caption{
\label{label}
We will utilize an abstract simplex in order  to {\it label} the hyperbolic tetrahedra which arise in this paper. 
To do so, view the simplex as the subsets of $\{1,2,3,4\}$, and 
identify each subset in the form $\{i\}$ with a tetrahedral face. The edge determined by $\{i\}$ and $\{j\}$ will be identified with the subset   $\{i,j\}$, or $\{ij\}$ for short.   Similarly, the vertex  determined by $\{i\}$,$\{j\}$,  and $\{k\}$ will be identified with the with the subset   $\{i,j,k\}$,  or $\{ijk\}$ for short.   
We will orient our simplex  and let $\left<\{1\},\{2\},\{3\},\{4\}\right>$ denote the positive orientation class.
 As such, for  any finite or ideal hyperbolic  tetrahedron, to  each $\{ij\}$  will may assign a
dihedral angle, $\theta(\{ij\})$.  
We will  call 
 $e^{i\theta(\{ij\})}$ the
{\it circulant } associated to $\{ij\}$ and 
we will  call  $e^{2 i \theta(\{ij\})}$ the {\it clinant  } associated to $\{ij\}$.    
We have labeled the {\it tetrahedral clinants},  and will denote them as 
 $( a,b,c,d,e,f )$, while the circulants will be denoted as 
  $\vek{c}=(\sr{A},\sr{ B},\sr{C},\sr{D},\sr{E},\sr{F})$ and the dihedral angles
   will be denoted as   $\vek{a}=(\di{A},\di{ B},\di{C},\di{D},\di{E},\di{F})$. 
  }
\end{figure}

Our scissors classes will  always  be closely related to  a collection of oriented ideal tetrahedra. In this section, we 
 fix some notation related to describing the  ideal tetrahedron and  review several useful facts concerning the   ideal tetrahedron.

Recall that an unoriented ideal tetrahedron is 
the convex hull of 4 ideal points in hyperbolic space  $H^3$.
 In order to describe this tetrahedron's orientation, and put coordinates on the collection of all such tetrahedra, it is useful to label this  unoriented ideal tetrahedron.  To do so we may label the ideal tetrahedron's four faces  as in figure \ref{label}.  For an ideal tetrahedron, this turns out to be more labeling than necessary, namely we can form coordinates which will only depend on the orientation class of the $\left<\{1\},\{2\},\{3\},\{4\}\right>$  labeling and the specification of an edge. 
 We can do this via the ideal tetrahedron's {\it complex coordinate}, $z(\{ij\})$, associated to  the  specified  $\{ij\}$ edge. To compute $z(\{ij\})$,  note that the sphere at infinity is naturally the Riemann Sphere, and  $z(\{ij\})$ can be computed as  the cross ratio of our ideal tetrahedron's four ideal points as
 \begin{eqnarray*}
z(\{ij\}) &   = & (\{ijl\},\{ijk\} ;  \{jkl\},\{ikl\}),
\end{eqnarray*}
where $\left<\{i\},\{j\},\{k\},\{l\}\right>$ is in $\left<\{1\},\{2\},\{3\},\{4\}\right>$ positive orientation class and
where  we have labeled our four ideal points  as in figure \ref{label}.

If $z(\{14\})$  is denoted as $z$, then the edge choice dependency is determined as follows
\begin{eqnarray}\label{c0}
z(\{14\})=z(\{23\})  = z,  
\end{eqnarray} 
 \begin{eqnarray}\label{c1}
z(\{12\})=z(\{13\}) = (z-1)/z, 
\end{eqnarray} 
\begin{eqnarray}\label{c2}
z(\{13\})=z(\{24\}) = 1/(1-z).  
\end{eqnarray}

Such an ideal tetrahedron has a natural orientation which can be easily determined by looking at the sign of 
any $\Im(z(\{ij\})$; if $\Im(z(\{ij\})$ is positive then the  ideal tetrahedron is positively oriented and if $\Im(z(\{ij\})$  is negative then the  ideal tetrahedron is negatively oriented.  If $\Im(z(\{ij\})$ is $0$, then  the tetrahedron is degenerate, see the note at the end of this section.  These concepts are independent of the choice of $\{ij\}$. Any permutation of the $\{i\}$ labels associated to
 $\left<\{1\},\{2\},\{3\},\{4\}\right>$'s negative orientation class,  negates the orientation.
Using our above $z$, after reversing the orientation class, we  have that 
\begin{eqnarray}\label{cross2}
 z(\{14\})=z(\{23\})  = 1/z, 
 \end{eqnarray}
\begin{eqnarray}\label{cross}
  z(\{12\})=z(\{34\}) =(1-z), 
\end{eqnarray}
\begin{eqnarray}\label{cross1} 
 z(\{13\})=z(\{24\}) = (z-1)/z.
\end{eqnarray}

In order to compute  the orientation-sensitive volume of the ideal tetrahedron,  we introduce the dilogarithm function 
\[\scr{L}_{2}(z)  =  - \int_{0}^{z} \frac{\log(1-s)}{s} ds, \]
which can be viewed as an analytic function with a branch cut along  $[1,\infty]$;
and  the Bloch-Wigner dilogarithm function 
\[ \scr{B}(z) = \Im(\scr{L}_2(z)) + \arg(1-z) \log|z|. \]
 The volume is given by
\begin{eqnarray}{\label{iv1}}
 \vol([z(\{ij\}]) &  =  & \scr{B}(z(\{ij\}),
   \end{eqnarray}  
 as  derived by S. Bloch and D. Wigner and presented by Milnor in  \cite{Mi2}.

We can express the complex coordinate of an ideal tetrahedron in terms of the ideal tetrahedron's clinants  via the 
the following observation 
\begin{eqnarray}{\label{fund}}
z(\{ij\}) & = & 
\frac{1-e^{-2i \theta(\{ik\}) } }{1-e^{2i \theta(\{il\}) } },
\end{eqnarray}
where $\left<\{i\},\{j\},\{k\},\{l\}\right>$ is in $\left<\{1\},\{2\},\{3\},\{4\}\right>$ positive orientation class. 
Using the notation   from  figure \ref{label}, a labeled oriented  ideal tetrahedron is equivalent to a list of clinants $(a,b,c,a,b,c)$ where $abc=1$, which will be denoted as $(a,b,c)_{it }$.  
Once again, this coordinate is determined by the orientation class of the facial labeling  and the  specified edge, which  corresponds to 
 $(a,b,c)_{it}$'s first coordinate.  The edge changes in equations (\ref{c0})-(\ref{c2}) correspond to even permutations of the  
 $(a,b,c)_{it}$ coordinates, while orientation class changes in equations (\ref{cross2})-(\ref{cross1})  correspond to 
the  odd permutations together with the conjugation of all the clinants. 
Given any  orientation reversing hyperbolic  isometry $I$,  
\begin{eqnarray}{\label{oreniso}}
 I((a,b,c)_{it }) & =  & (\bar{a},\bar{b},\bar{c})_{it } 
\end{eqnarray}
 and, hence, by equation (\ref{mirror})  
\begin{eqnarray}{\label{neg}}
 [(\bar{a},\bar{b},\bar{c})_{it }]  & =  &  -[(a,b,c)_{it }].
 \end{eqnarray}
Notice this implies that { all} permutations of the  $(a,b,c)_{it}$ coordinates preserve $[(a,b,c)_{it }]$, though the odd permutations require  equation (\ref{mirror}) and correspond to the original tetrahedron's positively oriented mirror image.
In terms of the $z$ coordinates,  $I(z(\{14\})) = \bar{z}$,  while  $1/\bar{z}$ corresponds to the tetrahedron's  positively oriented mirror image.

A particularly interesting special case of the ideal tetrahedron is the
\emph{isosceles} ideal tetrahedron
\[ (d^2,-1/d,-1/d)_{it }.\] 
Equivalently,
an ideal tetrahedron is isosceles
if its complex coordinate  is unit sized. 
In fact, for  $(d^2,-1/d,-1/d)_{it }$  we have  $z(\{14\})=d$, hence,  from equation (\ref{iv1}),
  \[ \vol([(d^2,-1/d,-1/d)_{it }]) = \Im( \scr{L}_2(d)). \]
It is geometrically straight-forward (see Milnor \cite{Mi1}) to prove that
 \begin{eqnarray}{\label{mis}}
2  [(a,b,c)_{it}]  = [(a^2,-1/a,-1/a)_{it },(b^2,-1/b,-1/b)_{it },(c^2,-1/c,-1/c)_{it }],  
  \end{eqnarray}    
hence 
 \begin{eqnarray}\label{iv2}
 \vol(2 [ (a,b,c)_{it }])& =  & \Im(\scr{L}_2(a)+    \scr{L}_2(b)+    \scr{L}_2(c)). 
  \end{eqnarray}    

{\bf Comment:} Letting $\scr{D} = \frac{arg(d)}{2}$, we have that
 \[  \frac{1}{2}  \Im(\scr{L}_2(d)  ) = \Lambda(\scr{D} ),  \]  
 where  $\Lambda(x)$ is the Lobachevsky function.
 Hence, as observed by Milnor (see  \cite{Mi1}),
 we may express equation (\ref{iv2}) as   
 \begin{eqnarray*}  \vol([ (a,b,c)_{it }])& =  & \Lambda(\scr{A}) +   \Lambda(\scr{B}) +   \Lambda(\scr{C}).  \end{eqnarray*}

{\bf Note:} If we compactify the space of ideal tetrahedra there are three types of degenerate tetrahedra, the {\it flattened} ideal tetrahedra, where $z(\{ij\}) $ is real and not in $\{0,1,\infty\}$,  and the {\it stretched} ideal tetrahedra, corresponding to a permutation of $(1,d,1/d)_{it }$ where $d \neq 1$, and the stretched and flattened tetrahedra where $d=1$ and $z(\{ij\})  \in \{0,1,\infty\}$.    The space of all oriented ideal tetrahedra with a specified edge 
can be described by the blowing up of  $S^1 \times S^1$ with center the point $(1,1)$.  In other words the real analytic  variety  formed by taking  $S^1 \times S^1$  and replacing $(1,1)$ with the set of lines through $(1,1)$.
To accomplish this we view    $(a,b) \in S^1 \times S^1$  as a pair of clinants.  If $a \neq 1$ and $b \neq 1$, then  $ (a,b,c)_{it} =(a,b,1/(ab))_{it}$.  When $(a,b)=(1,1)$ the slopes of the lines through   $(1,1)$  correspond to the real valued complex coordinates via equation (\ref{fund}).  


\end{subsubsection}

\begin{subsubsection}{ The Murakami-Yano formula \label{MY}}

Here we present the Murakami-Yano formula using the notation from sections \ref{siscon} and \ref{itet}. 
We will use the   tetrahedral circulants  as described in figure \ref{label}, and, throughout this subsection, and the next, we  let 
 $\vek{c}  =(\sr{A},\sr{ B},\sr{C},\sr{D},\sr{E},\sr{F})$ be the  tetrahedral circulants of a finite hyperbolic tetrahedron, which we shall denote as  $(\vek{c})_{tet}$.    Let $p=\sr{A}\sr{B}\sr{C}$,  $q =\sr{D}\sr{E}\sr{F}$,
 \[ \alpha(\vek{c})  =2 \left( \sr{A}\sr{D}+ \sr{B}\sr{E}+ \sr{C}\sr{F} + pq  + p \left( \sr{D}/\sr{A}+\sr{E}/\sr{B}+\sr{F}/\sr{C}+q/p  \right) \right),\]
\[ \beta(\vek{c})  = (\sr{D}/\sr{A}+\sr{E}/\sr{B}+\sr{F}/\sr{C}+\sr{A}/\sr{D}+\sr{B}/\sr{E}+\sr{C}/\sr{F})
-(\sr{A}\sr{D}+ \sr{B}\sr{E}+ \sr{C}\sr{F} +1/(\sr{A}\sr{D})+ 1/(\sr{B}\sr{E})+ 1/(\sr{C}\sr{F}) ),\]
and 
\[ \delta(\vek{c})  = |\alpha|^2-\beta^2. \]
In section \ref{clis}, we shall find that  $\delta$ is always positive, hence, using the positive square root,
  we may
 define the unit sized complex number   
\begin{eqnarray}\label{drho}
  \rho(\vek{c}) & =  & \frac{-\beta - i \sqrt{\delta}}{{\alpha}}. 
  \end{eqnarray}
Using $\rho$, we define 
the following list of isosceles ideal tetrahedra 
 \begin{eqnarray*} 
(\vek{c})_{\sr{M}\sr{Y}} &= & \left(\overline{\rho},-{\rho}\sr{C}\sr{D}\sr{E}, \overline{\rho \sr{B}\sr{C}\sr{E}\sr{F}},-{\rho}\sr{A}\sr{E}\sr{F}, \overline{\rho \sr{A}\sr{C}\sr{D}\sr{F}},  - {\rho}\sr{A}\sr{B}\sr{C},  \overline{\rho \sr{A}\sr{B}\sr{D}\sr{E}},-{\rho}\sr{B}\sr{D}\sr{F} \right).   
 \end{eqnarray*}
Murakami and Yano  derive the following formula (see \cite{Mu}):
  \begin{eqnarray}{\label{MYF}}
 \vol \left(4\left[(\vek{c})_{tet}\right]  \right) & = &  \vol \left(  \left[\right(\vek{c}\left)_{\sr{M}\sr{Y}} \right]+\left[\left(\bar{\vek{c}}\right)_{\sr{M}\sr{Y}}\right]\right) . 
\end{eqnarray}
 
We give an alternate proof of the Murakami-Yano formula by proving the following theorem. 
\begin{theorem}{\label{fundsis}}
Generically 
\[ 4\left[(\vek{c})_{tet}\right]  =     \left[\right(\vek{c}\left)_{\sr{M}\sr{Y}} \right]+\left[\left(\bar{\vek{c}}\right)_{\sr{M}\sr{Y}}\right].  \]
\end{theorem}
{\bf Sketch of Proof:} 
Here we describe how to geometrically realize this scissors congruence. 
That we arrive at  the correct angles, and that all the steps in this construction are well defined,  is verified  in section \ref{con}.    Throughout this sketch, terminology is used which hopefully will be clear to the reader from the indicated figures.

The first step in this construction is   summed up in figure \ref{cong} and its caption, 
which is a geometric summary of part \ref{p3} of theorem \ref{tetl} from section \ref{simp}.
In figure \ref{cong},  we see that  the $2 [(\vek{c})_{tet}]$ in the first row is equivalent to the scissors class of 
  the  pair of octahedra in the bottom row. We will call 
 these octahedral  a pair of  {\it \opc}, a concept we develop carefully in  section \ref{clis}. 
  It is also convenient to name the regions   in the second row of figure \ref{cong}.  We will call these regions  {\it supertetrahedra}, as introduced in figure \ref{suptets}.
 Notice in   figure \ref{cong}, that we utilized the  cutting down  procedure witnessed in figure \ref{heart}   to {\it cut down} our  supertetrahedra in row 2  to the \op in row 3.  In describing this cutting down procedure, we utilize  the convex version of the supertetrahedron.     This ability to use the convex case as a template for our constructions will prove very useful, and, in section \ref{simp}, we develop the notion a {\it $C$-region} to carefully justify this technique.  In order to appreciate  the utility of not needing to explicitly work outside the convex case,  we invite the reader to attempt to explicitly perform the cut downs in going from the supertetrahedra in the second row  of figure \ref{cong} to the octahedra in the third row.  
 
 Once we have our \opc,   we can perform the  {\it puff-and-cut} developed in  figure \ref{puffncut} with respect to any choice of vertex and face at this specified  vertex.  Performing this operation simultaneously to our \op in row 3 of figure \ref{cong}, utilizes a pair of oppositely oriented ideal tetrahedra pairs.   Hence, this puff-and-cut does not affect the  scissors class of our \opc.   In fact, a puff-and-cut is independent  of whether we use our specified face or the face opposite to our specified  face  at our specified vertex.  Hence, using the shading in  figure \ref{fireoct}, we can index our puff-and-cuts  as $P_v^s$ or $P_v^u$, depending on whether we use the shaded or unshaded pair of faces at $v$.  If we apply 
 \begin{eqnarray}\label{os}
  \osa & = & P^s_{\{14\}} P^s_{\{13\}} P^u_{\{34\}} P^s_{\{14\}} P^s_{\{23\}} P^s_{\{34\}}  P^s_{\{12\}}
  \end{eqnarray}
  to our \op  in the final row of figure \ref{cong}, 
  then, as a scissors class,  we have produced   8 ideal tetrahedra.
  At this point, we double these 8 ideal tetrahedra,  and, use equation (\ref{mis}) to express these doubled ideal tetrahedra as 24 isosceles ideal tetrahedra.  We  find that 8 of these 24 isosceles ideal tetrahedra  occur as pairs consisting of an  isosceles ideal tetrahedron together with an  oppositely oriented copy of this isosceles ideal tetrahedron. Hence $4 [(\vek{c})_{tet}]$ is equal  to the remaining 16  isosceles ideal tetrahedra.
  In section \ref{con}, we compute the angles arising in these 16  isosceles ideal tetrahedra and find these tetrahedra  are    precisely $ \left[\right(\vek{c}\left)_{\sr{M}\sr{Y}} \right]+\left[\right(\bar{\vek{c}}\left)_{\sr{M}\sr{Y}}\right]$, our needed scissors class.

 \qed

 \begin{figure}
\center{\includegraphics{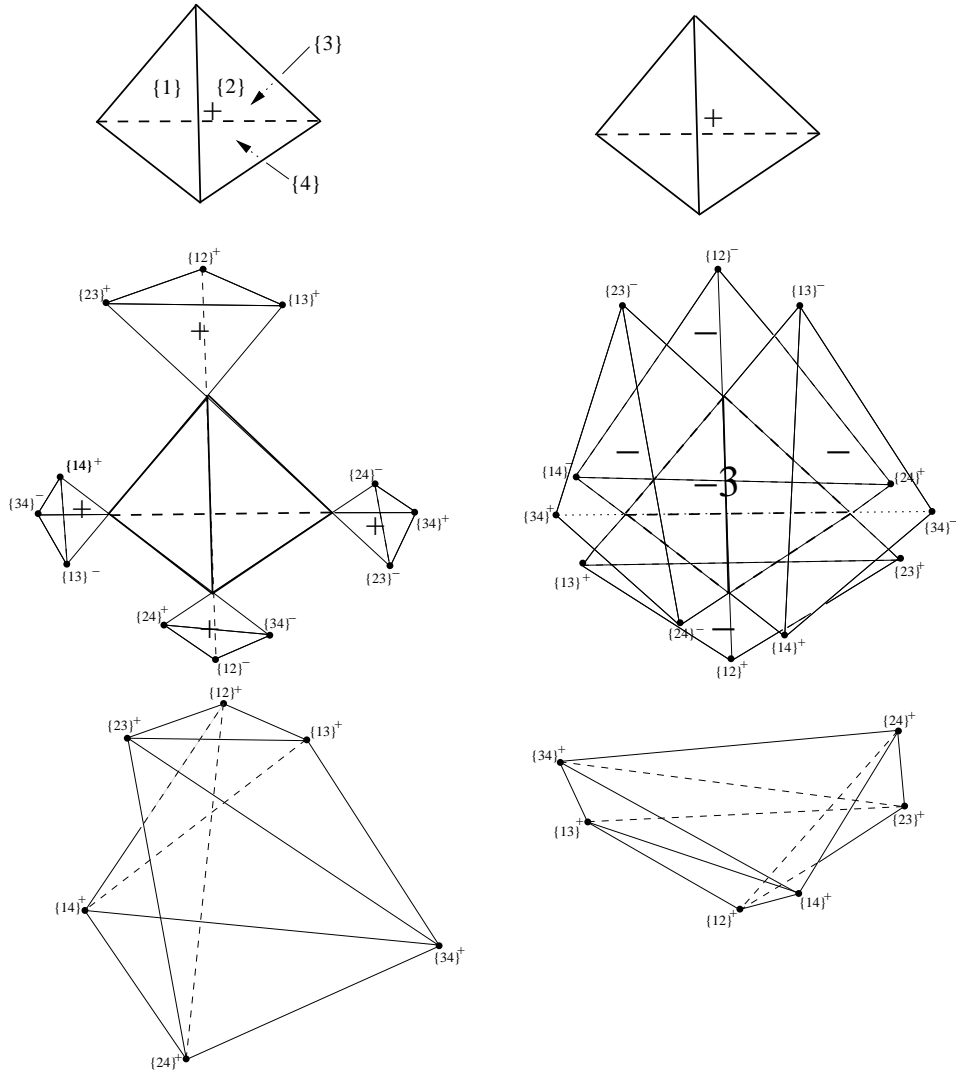}}
\caption{
\label{cong}
We start in the top row with two copies of our finite  tetrahedron.  
We then adjoin to the left hand copy the four indicated tetrahedra where the 
 bold faced vertices are ideal. 
To the right hand copy we remove the four indicated  tetrahedra.  The four adjoined and the four removed tetrahedra cancel out as a scissors  classes by  equation (\ref{mirror}).
The resulting regions in the second row  
can be cut down to the octahedra with the $\{ij\}^+$ vertices by removing the  6 ideal tetrahedra corresponding to the $\{ij\}^-$ vertices.   The combinatorics of our removal process is seen in figure \ref{heart}, where we  find that the ideal tetrahedra being removed  cancel in pairs.
Hence, we have realized  $2 [(\vek{c})_{tet}]$
as the scissors class of  the  pair of octahedra in the bottom row. 
}
\end{figure}

{\bf Comment 1:}   The term `generically' in the statement of
theorem \ref{fundsis} refers to the fact that 
theorem \ref{fundsis} is only proved for the finite hyperbolic tetrahedra in an open dense set 
 of the space of finite hyperbolic  tetrahedra.    In section \ref{dis} we describe the exact restrictions on our tetrahedra.  Despite this caveat, by continuity of volume,  theorem \ref{fundsis}   implies that equation (\ref{MYF}) will hold for all tetrahedra.

{\bf Comment 2:} Murakami and Yano
expressed the needed 
$\rho(\vek{c})$ and  $\rho(\bar{\vek{c}})$ using the roots of 
 the quadratic  polynomial
\begin{eqnarray*} \label{h1}
 h(\vek{c},z)
 & = & {\alpha} z^2 + 2 \beta z  + \bar{\alpha}  
 \end{eqnarray*}
or, more precisely,  the roots of
 $  \frac{\sr{A}\sr{B}\sr{C}\sr{D}\sr{E}\sr{F} }{2}h(\vek{c},z)$,
  which can be written as 
\begin{eqnarray*} \label{h2}
  \frac{\sr{A}\sr{B}\sr{C}\sr{D}\sr{E}\sr{F} }{2} h(\vek{c},z) &  =  & - 1/z  \left(  (1-z)(1-\sr{A}\sr{B}\sr{D}\sr{E}z)(1-\sr{A}\sr{C}\sr{D}\sr{F}z)(1- \sr{B}\sr{C}\sr{E}\sr{F}z) \right. \\ 
 & & \left.  -(1+ \sr{A}\sr{B}\sr{C}z)(1+\sr{A}\sr{E}\sr{F}z)(1+ \sr{B}\sr{D}\sr{F}z)(1+ \sr{C}\sr{D}\sr{E}z)  \right).
   \end{eqnarray*}
Notice that the two roots of $h(\vek{c},z)$ are $\rho(\vek{c})$ and  $\overline{\rho\left(\bar{\vek{c}}\right)}$.

\end{subsubsection}

\begin{subsubsection}{Alternate Views of Hyperbolic Volume \label{alt}}

One of the achievements of the Murakami-Yano formula was that it reduced the algebraic difficulty of volume computation to a quadratic equation.
Namely, the Murakami-Yano formula reduces the algebraic difficulty of computing volume to choosing  the appropriate square root of $\delta(\vek{c})$ in equation (\ref{drho}).  We now describe how to  analytically factor this  $\sqrt{\delta}$ out of the volume equation altogether, further revealing the algebraic  and analytic simplicity  of hyperbolic volume.   To do so, let $\gamma_i$    be  the $i^{th}$  component of the following vector
 \[ \vec{\gamma}=[1,-\sr{C}\sr{D}\sr{E}, \sr{B}\sr{C}\sr{E}\sr{F},  -\sr{A}\sr{E}\sr{F},\sr{A}\sr{C}\sr{D}\sr{F},-\sr{A}\sr{B}\sr{C},\sr{A}\sr{B}\sr{D}\sr{E}, -\sr{B}\sr{D}\sr{F}] . \]
 and
 \[ c_i = \frac{-i \gamma_i}{{\alpha}+\gamma_i \beta}. \] 
We find there is  an analytic function $\scr{F}$ such that  
\begin{eqnarray} \label{nicev}
\vol([(\vek{c})_{tet}]) &  = &  \Im\left( \frac{\sqrt{\delta}}{2}  \sum_{i=1}^{8} (-1)^{i}
c_i \scr{F}\left(\delta c_i^2
 \right) \right). 
\end{eqnarray}
 To describe $\scr{F}$, let
\begin{eqnarray*} 
 \scr{H}(w) & = & \scr{L}_2\left( \frac{2w}{1+w}\right)+  \frac{1}{4} \left( \log\left( \frac{1-w}{1+w}\right) \right)^2 
 \end{eqnarray*}
with branch cuts  $[-\infty,1] \bigcup [1,\infty]$ and $H(0)=0$.
In section \ref{dhnice}, we will see that 
\begin{eqnarray} \label{hnice}
 \vol(2 [(\vek{c})_{tet}])  &  =  &  \Im\left( \sum_{i=1}^{8} (-1)^{i}  \scr{H} ( c_i  \sqrt{\delta} ) \right). 
\end{eqnarray}
In section \ref{FSymd}, we will check that $H(w)$ is odd, which immediately  implies   equation (\ref{nicev}).
\end{subsubsection}

In section,  \ref{revis} 
several other views of the formula for hyperbolic volume are explored.

\end{subsection}

 \begin{subsection}{Scissors Congruences of the Tetrahedron  }
 
 At the heart of the proof of theorem \ref{fundsis} is the construction of the transformation  $\osa$ from equation (\ref{os}). 
 The puff-and-cuts used to build $\osa$ generate  a group of scissors congruences, which we will discuss in this section.
To make good sense out this group,  we first develop the notion of a   
  {\it   generalized hyperbolic  tetrahedron}.

 \begin{subsubsection}{The Generalized Hyperbolic  Tetrahedron  \label{com}}

Generalized hyperbolic  tetrahedra are  oriented, labeled subsets of  hyperbolic space, that include,  as a special case, the  finite hyperbolic  tetrahedra.  The generalized hyperbolic  tetrahedra should be viewed as the natural analytic continuation of  the space of  finite hyperbolic   tetrahedra.    To make these  tetrahedra less mysterious we have the following result.

 \begin{theorem}\label{tetplane}
Up to orientation preserving isometry, every collection of four, distinct, labeled, pairwise intersecting planes in  $H^3$ corresponds to a unique  generalized hyperbolic tetrahedron.  
 \end{theorem}
 {\bf Sketch of Proof:}  
 This theorem is  proved in section \ref{simp}, where the generalized hyperbolic  tetrahedron is carefully developed.   Throughout this sketch, various terminology is used which hopefully will be clear to the reader from the indicated figures. 
 
Notice that the labeling scheme from figure \ref{label} can be used to label   
any collection of four, distinct, labeled, pairwise intersecting planes in  $H^3$.
 As such, we call each plane $\{i\}$ a  {\it tetrahedral face}, each  geodesic $\{ ij \}$  when $i \neq j$ a {\it  tetrahedral edge} and   each $ \{ijk\}$ a {\it    tetrahedral vertex}, when   $ \{ijk\}$ is nonempty and $i$,$j$,and $k$ are distinct.  
   We will find that  a bit more labeling is sometimes desirable, and we will say that  such a labeled collection of planes  is  {\it decorated} if 
   each geodesic  $\{ij\}$ has been assigned an orientation. 
   This is equivalent to labeling the two end points of $\{ij\}$ at infinity as  $\{ij\}^+$ and  $\{ij\}^-$.  
       
Our first step of this proof will be  to construct  the
   general supertetrahedron, a concept introduced in the sketch of theorem \ref{fundsis}'s proof. 
To do so,  take any   four planes  as described in the statement of the above theorem 
 and decorate them.
 A supertetrahedron can be constructed by attaching to the $\{ij\}^{\pm}$  ideal vertices  the  ideal tetrahedra witnessed in figure \ref{heart}.
 A couple of examples of the sorts of regions that can result from this procedure are seen in the second row of figure \ref{suptets}.

From our supertetrahedron we can construct our needed  generalized hyperbolic  tetrahedron.
 We can motivate this  construction    by 
 attempting to reverse, for any supertetrahedron, the  process of  turning the  finite tetrahedra in row 1 of figure \ref{cong}  into the supertetrahedra in row 2.  The four   tetrahedra utilized in figure \ref{cong} to go from row 1 to 2 are special cases of {\it half prisms}, see the second row of figure \ref{prisms}.
  In general, at each  `vertex' of a supertetrahedron there is well defined prism, see figure \ref{prisms}.
Starting with a supertetrahedron, we can remove the top halves of these four prisms and call the resulting region a 
{\it generalized hyperbolic  tetrahedron}, see figure \ref{tets}.  As in figure \ref{cong},  we see that the finite tetrahedron is a special case of this construction.

At this point in the proof,  we have associated a  generalized hyperbolic tetrahedron to every collection of collection of four, distinct, decorated, pairwise intersecting planes in  $H^3$. To  finish this proof,  we need to demonstrate that the region determined by this construction is independent of  our the edge orientations provided by our decoration.  This fact is  discussed in figure \ref{boundary}, and
 figure \ref{boundary}'s caption  completes our sketch. 

\qed

\begin{figure}
\center{\includegraphics{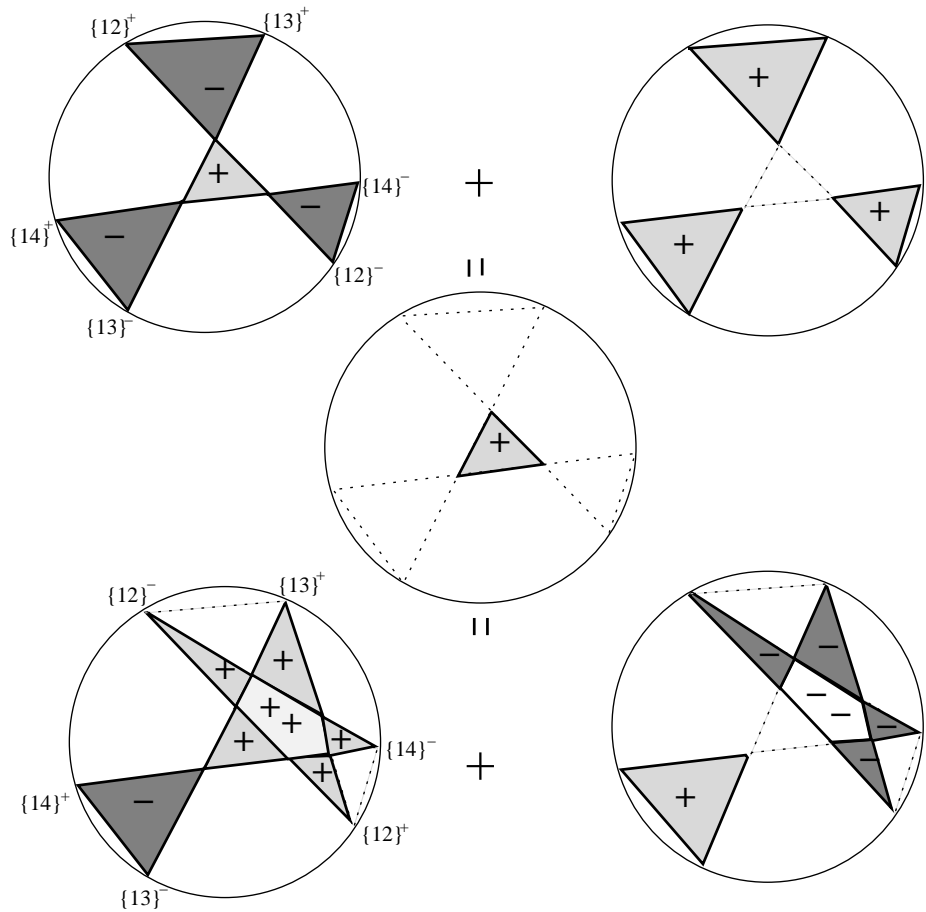}}
\caption{
\label{boundary}
Here we describe why the construction of our  generalized hyperbolic  tetrahedron is independent of  our $\pm$ sign choices. We will call the process of changing  such a choice the process of {\it flipping an edge}.
The flipping of an edge simultaneously modifies the supertetrahedron and the four  prisms involved in the tetrahedron's construction.  In the top row of his figure, we see  what the supertetrahedron and three prisms  on the left hand side of the second row in figure \ref{cong}, look like in the $\{1\}$ plane before flipping the $\{12\}$ edge. In the bottom row of this figure, we see  what the supertetrahedron and three prisms  look like after  flipping the $\{12\}$ edge.
Notice that in both cases, our supertetrahedron and half prisms combine to together to form the tetrahedral face in the center of this figure.      The fact that this is true in all the $\{i\}$ planes implies that flips preserve the needed tetrahedron (since two such regions that share the same boundary will share the same interior). 
These ideas will be carefully developed in section \ref{simp}. 
}
\end{figure}

 Theorem \ref{fundsis} holds, generically, for all generalized hyperbolic tetrahedra. 
From this point on,  we will refer to a labeled generalized hyperbolic tetrahedron as a {\it hyperbolic tetrahedron}, or simply a {\it tetrahedron}. When we want to emphasize that a hyperbolic  tetrahedron is finite, ideal or decorated we will say so.
We will find that every tetrahedron can be assigned angle data as described in figure \ref{label}.
 We will let {\it  tetrahedral angle data}  refer to all  six-tuples  of real numbers modulo $2 \pi$ that correspond to a 
 tetrahedron's associated dihedral angles, denoted as  $\vek{a}=(\di{A}, \di{B},\di{C},\di{D},\di{E},\di{F})$. 
  This association will not be unique.  Generically, there are   $2^8$ different collections of  tetrahedral angle data that correspond to the same  tetrahedron.  
 
 {\bf Comment:}   
 Using the planes associated to our tetrahedron, we may  
 decorate a tetrahedron by decorating its associated planes, as defined in the sketch of theorem \ref{tetplane}'s proof. 
 The space of  decorated tetrahedra form a $2^6$ fold cover of the space of   tetrahedra.
  Decorating a tetrahedron resolves most of the ambiguity in associating  tetrahedral angle data to its
edges. Namely, there are only  $4$ different collections of  tetrahedral angle data that correspond to the same  decorated tetrahedron.  
They are all related by the  copy of  the Klein 4 group  generated by adding $\pi$ to all but an opposite pair of dihedral angles, or rather the negation of all but an opposite pair of circulants.  We shall call this group $\trivial$.  Geometrically, 
$\trivial$ plays no and in order  to produce honest coordinates on the space of decorated tetrahedra will require modding out by $\trivial$'s action.
Such coordinates are useful for a variety of technical reasons, and in section \ref{algsym} we introduce them as the {\it balanced coordinates}.

\end{subsubsection}
 
 \begin{subsubsection}{Volume Preserving Transformations \label{vsym}}

We will say that $\gen$, a transformation acting on  the space of tetrahedral angle data, as introduced in   section \ref{com}, is a {\it volume-preserving transformation}  if the volume of the tetrahedron with    tetrahedral angle data $\vek{a}$ is the same as the volume of the tetrahedron with  tetrahedral angle data $\gen \cdot \vek{a}$.
One of the interesting corollaries of the Murakami-Yano formula is that 
the formula  exposes the fact that the elements of the Regge group are
hyperbolic volume-preserving transformations.  The Regge
group was discovered originally as an order-144 group that preserves the
classical $6j$-symbol.  
 In order to describe  the Regge symmetries,  let 
 \[\di{S}  =   \frac{\di{B} + \di{C}+ \di{E}+ \di{F}}{2} .  \] 
 The  Regge Group, viewed as acting on tetrahedral angle data, consist of the group generated by the tetrahedral symmetries together with the transformation $g_r$ satisfying 
   \[ \regg \cdot \vek{a}= (\di{A},\di{S}- \di{B},\di{S}-\di{C},\di{D},\di{S}-\di{E},\di{S}-\di{F}).  \] 
  Regge and Pozano conjectured that from the
$6j$-symbols one can reconstruct the volume of a Euclidean polyhedron,
which  would imply that the Regge symmetries are, in fact, 
Euclidean volume-preserving transformations.    Roberts resolved the Regge-Pozano conjecture and also observed that
each of these symmetries is realized by a  scissors congruence (see \cite{Ro}).  In
particular, this allowed Roberts to construct 12 distinct tetrahedra scissors
congruent to a fixed Euclidean tetrahedron (generically).  To accomplish this
Roberts utilizes the Dehn sufficiency theorem in Euclidean space
(that volume and the Dehn invariant determine a scissors class, see Dupont and Sah \cite{Sa3}). In
hyperbolic space, Dehn sufficiency remains  one of the fundamental unsolved
conjectures concerning the nature of  hyperbolic volume (once again, see Dupont and Sah \cite{Sa3}). 
If the Hyperbolic Dehn Sufficiency Conjecture is true, then, from
the Murakami-Yano formula, the Regge symmetries would be realized by
scissors congruences.   
However, without a resolution to the  Dehn
sufficiency conjecture, proving that these symmetries are
realized by scissors congruences will require other tools. 
By an explicit construction, Mohanty  has proved that the Regge  symmetries are scissors congruences (see   \cite{Mo2} and the comment at the end of this section).    As a
scholium to our methods,  we find that the Regge symmetries are
just the beginning. In section \ref{group}, we  construct the  group described in the following result.

\begin{theorem}\label{res2}
There is  an order-$23040$ group of  scissors-class-preserving transformations  acting properly on the  space of decorated tetrahedra.   We will call this group $\wacky$.
\end{theorem}
 {\bf Sketch of Proof:}
In section \ref{group},  these constructions, and the reasoning behind them, is developed carefully. 

First, we describe how to generate  $\wacky$ and why $\wacky$'s elements induce scissors congruences.
In section \ref{algsym}, we give a careful algebraic description of $\wacky$.  In particular, we find that $\wacky$  is isomorphic to  $D_6$, where $D_6$ is a well known  order-$23040$ reflection group.    

Let $\wacky$ be the group generated by the  puff-and-cuts from figure \ref{puffncut}.
We will view this group as acting on the set of all \opc.  
\Op will have conjugate dihedral clinants, and, hence, the ideal tetrahedra involved in performing the puff-and-cut will pairwise cancel out, just as they did with respect the \op arising in sketch of theorem \ref{fundsis}'s proof. 
In order to  see how to apply  $\gen \in G$ to a decorated  tetrahedron, first apply the construction in figure \ref{cong}, to realize the double of our  tetrahedron  as a pair of \opc.  Notice, going from the first row to the second row   in figure \ref{cong}, is well defined in general  since our decoration allows us to canonically associate a supertetrahedron to our tetrahedron, as in the sketch of theorem \ref{tetplane}'s proof.   To these \op apply  $\gen$.
 This results in new \opc.    We then invert the construction in figure \ref{cong}, which produces a new (decorated) 
 doubled tetrahedron.  By construction, this new doubled tetrahedron  is scissors congruent to its underlying \opc.  Hence, the original doubled tetrahedron is scissors congruent to the new doubled tetrahedron, and,  by equation (\ref{half}), the tetrahedra themselves are scissors congruent.

\qed

In terms of volume-preserving transformations, we have the following corollary. 

\begin{corollary}{\label{cor0}}
There is an order-$92160$  group of
hyperbolic volume-preserving transformations
acting properly on the space of tetrahedral angle data.  
We will call this group $\wack$.
\end{corollary} 
\proof
We can lift   $\wacky$'s action to a proper action 
 on the   tetrahedral angle data corresponding to decorated tetrahedra.
 For example, 
  the $\osa$ from equation (\ref{os}) acts via
  \[\osa \cdot \vek{a}=  (-\di{A},-\di{S}, \di{E}+\di{F}-\di{S},-\di{D}, \di{B}+\di{E}- \di{S}, \di{B}+\di{F}-\di{S}).   \] 
This action determines  most of $\wack$. 
 However, as in the comment at the end of section \ref{com},
 the  space of  tetrahedral angle data corresponding to decorated tetrahedra forms a four fold cover of the 
space of  decorated tetrahedra,  hence we should throw in the deck group of this cover.  This is precisely  the group $\trivial$  from the comment at the end of section \ref{com}.   $\trivial$ commutes with our lifting of  $\wacky$'s action. 
 Hence all the elements of $\wack = \wacky \times \trivial$ correspond to volume-preserving transformations.

\qed

From the sketch of theorem \ref{res2}'s proof, we see that   $\wacky$ acts on the space of decorated tetrahedra. However, $\wacky$ does not act (as a group) on the space of   tetrahedra. 
This is because the edge flips, from figure \ref{boundary}, generate a subgroup of $\wacky$   which preserves any tetrahedron but is not normal in $\wacky$.   This group generated by the edge flips is isomorphic to  $\left( \frac{\bf Z}{2 {\bf Z}} \right)^6$, 
called the   {\it shaded subgroup}, and plays a fundamental in the proof of the following corollary.


\begin{corollary}{\label{cor}}
There are  generically 30
(generalized) tetrahedra scissors congruent to a fixed  tetrahedron no pair of which are congruent to each other via an orientation preserving isometry.
\end{corollary} 
\proof
$\wacky$ contains many elements which preserve a  tetrahedron up to orientation preserving  isometry. For one  $\wacky$ contains the tetrahedral symmetries, see section \ref{sispre}. 
$\wacky$ also contains the shaded subgroup, which,  
 in terms of its action on the tetrahedral angle data,   corresponds to the negation of the individual coordinates.
 Let $P$ be the order-$768=12 \cdot 2^6$ subgroup group of  $\wacky$ generated by the edge  flips and the
orientation preserving tetrahedral symmetries.   The elements of $\wacky/P$ are our 30 candidate  tetrahedra. 
To see that these tetrahedra are distinct from  our original tetrahedron, we  will need  that  the cosines of 2 times the  dihedral angles  of each representative are well defined.  This is guaranteed by lemma \ref{baby23} and the fact that cosine is an even function.    To finish our proof, we look at representatives of each of these 30 cossets and note that none (but the identity) have these cosines of 2 times their  dihedral angles related to the  original tetrahedron by an orientation preserving tetrahedral symmetry.  To verify this requires looking at a representative of each of the 30 cosets of $\wacky/P$. In section \ref{sispre}, we write down  these cosets.

  \qed

  {\bf Comment:}  As discussed in section \ref{com},  among generalized  tetrahedra the tetrahedral angle data is not well defined.  However, among finite tetrahedra this concept is perfectly well defined. 
The problem is that $\wacky$ will take finite tetrahedra to non-finite tetrahedra, and force us   to give up this notion. 
In fact, if we start with a finite tetrahedron, then only 12 of the 30 tetrahedra from corollary \ref{cor} are finite tetrahedra (see section \ref{dis}).
These twelve tetrahedra correspond precisely to the Regge scissors classes.
 Mohanty  produces a construction of the Regge symmetries, a construction where  the notion of dihedral angle remains well defined throughout the process  (see  \cite{Mo2}).  
 

\end{subsubsection}

\end{subsection}

\end{section}

\begin{section}{Constructions \label{con}}

  In this section, we carefully describe the geometric  constructions necessary to prove  theorems
\ref{fundsis}, \ref{tetplane}, and \ref {res2}. 
In our first section, we  develop the notion of a generalized hyperbolic tetrahedron,  as discussed in section \ref{com}.

\begin{subsection}{The Generalized Hyperbolic Tetrahedron \label{simp}}

Here we will construct the  generalized  hyperbolic tetrahedron  
needed to prove theorem \ref{fundsis}.  To do so we utilize the concept of a $C$-region.
A $C$-region gives us a way of building non-convex regions of ${ \bf H}^3$  utilizing a convex  ideal  polyhedron, $C$, as a template.  These $C$-regions are technically very convenient,  and, by utilizing them, one can construct an open dense subset of all hyperbolic tetrahedra, see comment 2 at the end of this section.   
Some of the results in this section  are most naturally proved 
by induction from the one dimensional case, hence we will also discuss 
these notions in   ${ \bf H}^1$ and  ${ \bf H}^2$.   In fact,  all the results in this section 
 have analogs in ${ \bf H}^n$.

{\bf $C$-Regions, the idea:} Intuitively,  the $C$-region  is simple to describe.  Realize $C$ as a chain of ideal tetrahedra.  Viewing this chain as an abstract chain, a $C$-region is another realization of this chain with ideal tetrahedra, such that, if 
a collection of $C$'s vertices all lived on a  face of $C$, then the vertices are still cohyperplanar in the $C$-region.
Examples of $C$-regions are given in figures  \ref{prisms} and \ref{suptets}, and are always denoted via a pair $|C,R|$, where $C$ is the convex template and $R$ denotes a realization of the template.  The geometric region determined by $|C,R|$ (denoted as $||C,R||$),  the scissors class  determined by $|C,R|$  (denoted as $[C,R]$)  and $|C,R|$'s dihedral clinants of are all well defined and independent of the chain used to realize $C$.   This intuitive picture, along with these facts, should be enough to understand the constructions in this chapter, but here is a formal definition.

{\bf $C$-Regions, the definition:} We will use the Klein model of
${ \bf H}^3$,  and, hence, may view 
any chain of hyperbolic tetrahedra with all finite
and/or ideal vertices  as a chain of Euclidean tetrahedra
in ${\bf E}^3$. 
Let $S$ be an  finite  oriented 3-dimensional
simplicial complex.  We let  a {\it realization}
$R$ of $S$ be an assignment  of a
point $R(p) \in {\bf E}^3$ for every vertex $p \in S$. A
pair $(S,R)$ will be called a {\it 
Euclidean chain} provided $R$ is injective. 
We may  simplicially extend
$R$ and continuously map our abstract complex into ${\bf E}^3$,
and,  assign an integer {\it label} to a  full measure set of points in  ${\bf E}^3$
corresponding to the local degree of this mapping.   This labeled set will  be called
the {\it region} determined by $(S,R)$.  If $R(p)$ is  in  the
closed unit ball,  then  $(S,R)$ may also be used to  represent a well
defined {\it Hyperbolic chain},  and the region determined by $(S,R)$  
can be thought of as a labeled subset of ${\bf H}^3$, which will be 
denoted as  $||S,R||$.   Let $[S,R]$ be the scissor class of the list of  tetrahedra 
  determined  by  $(S,R)$.   If all the vertices are ideal,  then we
call $(S,R)$ an {\it  ideal chain}. Utilizing the ideal tetrahedra used to form an ideal chain, we find  that 
  an  ideal chain has  a well defined scissors class and well defined dihedral circulants.
   
  If  $C$ is a convex ideal
polyhedron  we may  triangulate $C$ using its ideal vertices, and, hence, realize  $C$ 
 as $||S_C,R_C||$ for some ideal chain $(S_C,R_C)$.  Given an $(S_C,R)$, if  
any set of vertices that shared a top dimensional facet in $C$ 
are still cohyperplanar under $R$, then we will say that 
$R$ satisfies  $C$'s {\it facial constraints}.  We will let 
a {\it $C$-Region}, $|C,R|$, be an labeled subset of $H^3$
  which is $||S_C,R||$ for some $S_C$  where $R$ satisfies $C$'s facial
constraints and where  all of $C$'s   top dimensional facets determine 
distinct hyperplanes.  
  
\begin{lemma}\label{baby23} 
If $|C,R|=||S_C,R||$,
then any other chain, $\hat{S}_C$, used to realize $C$ we have that 
  $|C,R|=||\hat{S}_C,R||$.  
Furthermore, $|C,R|$'s scissors class,
 $[S_C,R]$, and dihedral  clinants are well defined, in other words, independent of the choice of  $S_C$. 
\end{lemma}
\proof
Our goal will be to take any two chains $S_C$ and   $\hat{S}_C$ and to show that 
 $[S_C,R] = [\hat{S}_C,R]$,   $||S_C,R]||= ||\hat{S}_C,R||$  and that the dihedral clinant agrees.
 
To begin, notice   $\hat{S}_C$ induces a triangulation of each of $C$'s faces.   To such a triangulation we'd like  to implement  a sequence {\it 22-moves}, where a 22-move  takes two triangles that meet in the diagonal of a quadrilateral and replaces them with the triangles forming the quadrilateral that share the quadrilateral's other diagonal.    It is simple to verify that every triangulation of a convex polygon utilizing only the polygon's vertices as vertices of the triangulation  is equivalent to every other such triangulation  via a sequence of such 
22-moves.
 In n dimensions we have the analogous 2n-move, for example the 23-move in figure \ref{two3} , which also has this property. 
 
 \begin{figure}
\center{\includegraphics{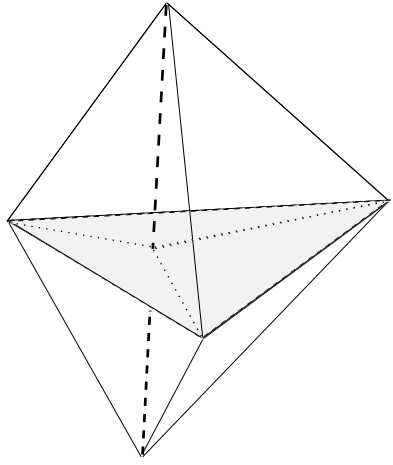}}
\caption{
\label{two3}
 Here we see how to view a pair of tetrahedra glued along the shaded face as the triple of 
 tetrahedra that share the bold edge, or conversely.   Such a view point change is an example of a 23-move.}
 \end{figure}

The next observation is that the  22-moves of  $\hat{S}_C$'s boundary can be implemented  by adding or removing tetrahedra to a face of $C$.  
Since $R$ satisfies $C$'s facial constraints these tetrahedra are degenerate and hence will  not change  $[\hat{S}_C,R]$,   $||\hat{S}_C,R||$ or any dihedral clinant.
This procedure may change a dihedral  circulant by $-1$, and hence the circulants will not be well defined.   Utilizing these degenerate tetrahedron, and the observation that all the facial triangulations of $C$ differ by 22-moves, assures us that we may, without loss of generality,   assume 
 that the triangulation of $\hat{S}_{C}$'s boundary agrees with
the triangulation of ${S}_{C}$'s boundary. 
 
 Since  $\hat{S}_{C}$'s and  ${S}_{C}$'s boundary agree,  
we may glue together  $S_C$ to $-\hat{S}_C$ to form a sphere. $R$
induces a mapping of  $S^3$ into $E^3$ which must have degree zero, hence,
$||S_C,R|| =||\hat{S}_C,R||$ from which $[S_C,R]= [\hat{S}_C,R]$ immediately follows. Since the sphere is boudaryless  the clinants around every edge must multiply to one,
hence the dihedral clinants of   $|S_C,R|$ and $|-\hat{S}_C,R|$ must be conjugate.
The clinants of  $|-\hat{S}_C,R|$  are conjugate of those of $|\hat{S}_C,R|$, hence, by equation (\ref{neg}),  the clinants of $|\hat{S}_C,R|$  and $|S_C,R|$ must agree.  

 An independent proof can be accomplished by  noticing that all triangulation of $C$ are equivalent via a sequence of the {\it 23-moves}  in  figure \ref{two3}.     In fact, the 23-move  forms the  primary relation in  $\scr{P}(H^3)$,  when $\scr{P}(H^3)$  is viewed as generated by ideal tetrahedra (see Dupont and Sah \cite{Sa1}).

\qed

 \begin{figure}
\center{\includegraphics{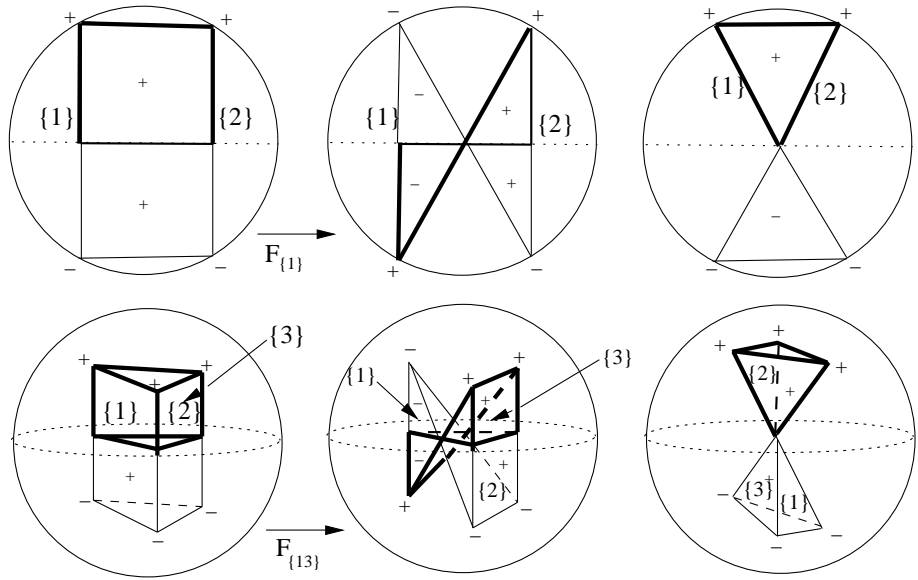}}
\caption{\label{prisms}
On the left, we have  2 and 3 dimensional  convex prisms, followed by
examples of convex prism-regions, which will be called {\it prisms}.  We
will denote such a prism as $|P,R|$. The convex prism in
$H^{3}$  can be defined by placing a hyperbolic triangle   in a hyperplane of $H^3$, and then taking the
convex hull of the union of the geodesics perpendicular to the hyperplane
that contains the vertices of this specified triangle. 
Via the labeling convention described in  section \ref{simp}, we can utilize the
pictured hyperplanes in order to label the vertices as
$\{ij\}^{\pm}$.  As such, the labeled   $\pm$ signs  are sufficient to decorate our vertices. 
  $||P,R||$ can be naturally cut in half, and 
each component of the interior of  these halves has either positive or negative vertices in its closure. 
Hence, every prism has a well defined {\it top} associated to 
its positive vertices and a well  and {\it bottom} associated to 
its negative vertices.   In the figure,
the top half is indicated with the bold face lines. 
In going from our left most to our middle prism we have 
demonstrated geometrically the notion of a flip, as introduced in section \ref{simp} .
 }
\end{figure}

\begin{figure}
\center{\includegraphics{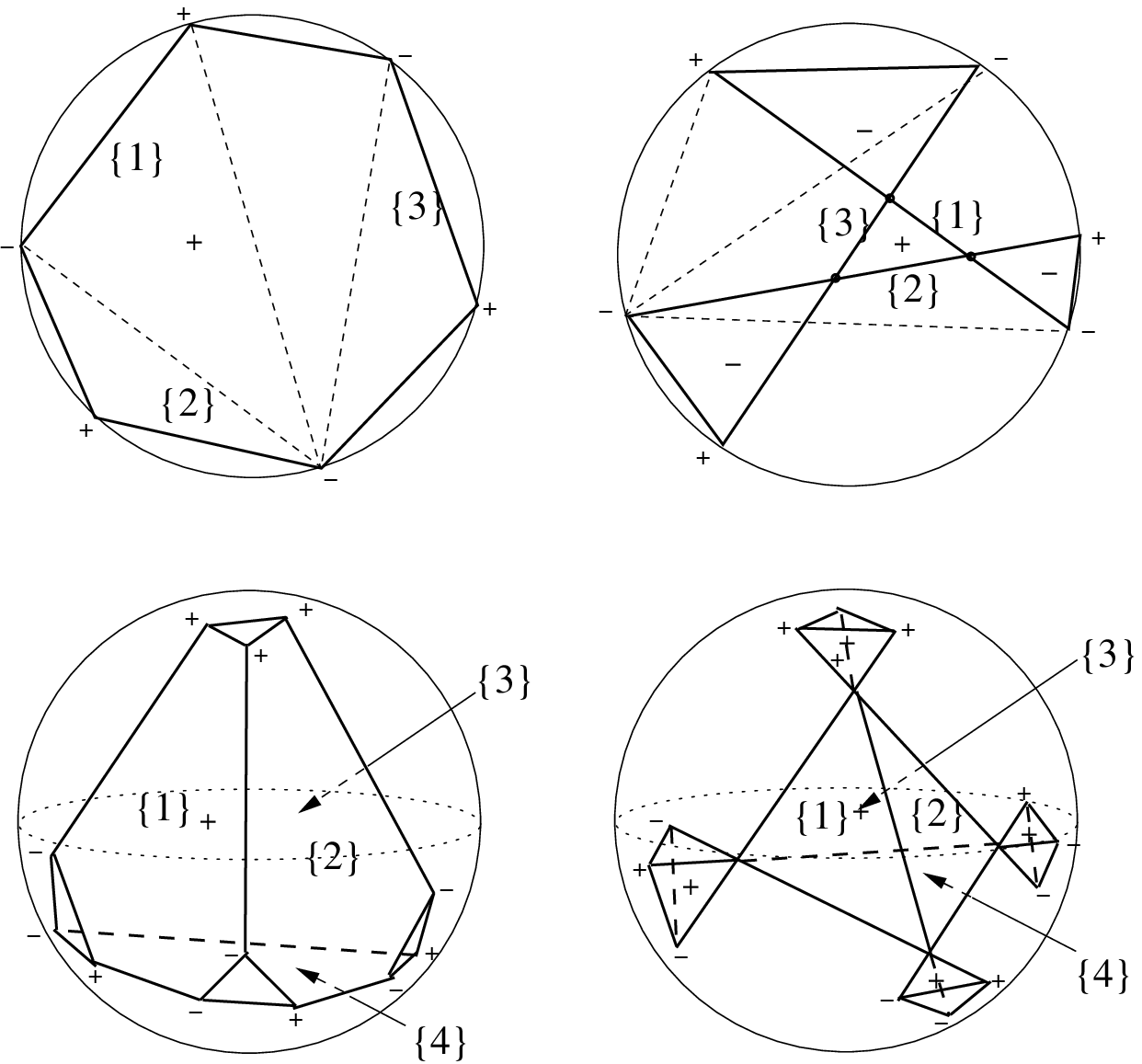}}
\caption{\label{suptets}
On the left, we have  2 and 3 dimensional  convex supersimplicies, followed by examples of convex supersimplex-regions, which
will be called {\it supersimplicies} and denoted $|Su,R|$. To construct a {\it supertetrahedron}, first 
take a Euclidean tetrahedron in the Klein model with all its vertices
hyperideal, and with the property that each of its edges intersects
$H^3$. The convex supertetrahedron can be defined as  the convex hull of the
geodesics  formed by  intersecting the Euclidean tetrahedron's edges with
$H^3$, as indicated in the figure.  Note, there are $4$ natural prisms in a supertetrahedron. 
We have utilized the $4$ labeled hyperplanes in order  to label our supertetrahedron, and
to each subset of  $3$ of these hyperplanes we  have an associated a prism. We will
let $|P,R_i|$ be the prism associated to the hyperplanes 
 $\{1,2,3,4\}-\{i\}$ with its top determined by the three vertices whose
convex hull is a facet of the convex supertetrahedron. 
}
\end{figure}

{\bf Hyperplane Notation:} In figures \ref{prisms} and \ref{suptets}, we see our most important examples of $C$-regions, the prism and the supertetrahedron. We use hyperplanes to indicate our method of labeling the vertices.  
We use the labeling scheme described in the sketch of theorem \ref{tetplane}'s proof. 
  We says $\{1,2,3,4\}$  is a {\it non-degenerate} collection  if each
  $\{ ij \}$ with $i \neq j$  is one dimensional and  if  each
    $ \{ijk\}$  fails to be  an ideal point.

    {\bf Flipping an Edge:} With such a labeling we have the notion of {\it flipping an edge}, described in figure \ref{boundary}. In terms of our $C$-region notation, flipping the edge      $\{ ij\}$   correspond to changing the roles of  $\{ ij\}^+$  and    $\{ ij\}^-$. In other words,  starting with   $|C,R|$ we form  $\left|C, F_{ij}R \right|$, where $R$ and $ F_{ij}R$ agree  on all the   vertices of $S_C$  accept    
    \[  F_{ij}R (\{ ij \}^{\pm}) =  R (\{ij \}^{-\pm}). \] 
    We define  $\Omega R$ be the new realization determined by flipping all of  $|C,R|$'s edges.
 
   \begin{lemma}{\label{prisl}}
The following facts are true about  the 3 dimensional prism described in figure \ref{prisms}.  
\begin{enumerate}
\item
A  prism is equivalent to a  collection $3$ decorated, non-degenerate  planes.
\item 
A prism can be divided in to a top and bottom half, which are scissors
 congruent .
\item
If we flip all the edges of a prism then 
the top half of the new prism is isometric to the
 top half of the  original prism, and via an orientation reversing isometry. 
\end{enumerate}
\end{lemma}
\proof
Given a prism $|P,R|$, since $R$ satisfies the needed planarity constraints,  we can construct the  planes needed in the first part of this lemma.  Conversely with these decorated planes  we have determined where  the   vertices of $C$ must go under $R$ which, by lemma \ref{baby23}, completely determines our needed prism. 

Given a prism the three planes it determines either intersect at a finite or hyperideal point  (since $R$ is injective the point cannot be ideal, see comment 2 at the end of this section).
Hence, with a hyperbolic isometry we can send this point to infinity if it is hyperideal or to the origin in the Klein model this point is finite. 
Then, via edge flips, we can arrange this configuration of planes to be qualitatively either the first or final prism in row two of figure \ref{prisms}.  In these two cases the prism's top and bottom are mirror images of each other, hence by equation (\ref{mirror}), the top and bottom scissors congruent.  Now we simply pick an ideal triangulation, and  perform the edge flips to see that  there is, up to orientation, only one other qualitatively   distinct case, that of  the middle prism in the second row of figure \ref{prisms}.  This case  also satisfies  this mirror image property as needed.   Notice that lemma \ref{baby23} assures us that we need not examine what takes place with regard to other triangulation choices. 
One can also  prove the second part by an induction from the 1 to the 2 to the 3 dimensional case and beyond.  However, in even dimensions, the bottom half of the finite vertexed prism will have the opposite orientation of the top. 

Having explicitly described all our prisms, the third part can simply be verified by  picking an  ideal triangulation,   and explicitly performing  the $\Omega$ transformation in the three qualitatively distinct cases.  This third part  is always true in odd dimensions.
 \qed

\begin{lemma}{\label{supl}}
The following facts are true about  the supertetrahedron described in figure \ref{suptets}.  
\begin{enumerate}
\item
A supertetrahedron is equivalent to a collection $4$ decorated, non-degenerate  planes.
\item 
Every supertetrahedron  determines $4$ prisms, and the region determined
by removing the top halves of these $4$ prisms is independent of
edge flips. 
\end{enumerate}
\end{lemma}
\proof
The first part is proved exactly as  the first part of lemma \ref{prisl} was proved.
That every supertetrahedron geometrically determines 4 prisms follows from the fact that our our supertetrahedron is equivalent to a  collection $4$ decorated planes, which provides  us with the four collections of $3$ decorated planes, needed to construct  our  prisms from lemma \ref{prisl}. 
 That the region formed 
by removing the top halves of these $4$ prisms is independent of
edge flips can be proved  much like the proof of lemma \ref{baby23}.
Namely,  we need only show that such a  region has the same boundary before and after an edge flip. There are two types of boundary components. The part contained in the planes determined by the supertetrahedron and the hyperbolic triangles in the waists of the prisms when  a triple of these planes intersect at a hyperideal point.   The  hyperbolic triangles in the second case
are easily seen to be edge flip invariant by simply performing the needed flip to each of the 2  qualitatively different  types of such prisms found in figure  \ref{prisms}.  
To prove that the parts of the boundary  contained in the planes determined by the supertetrahedron are preserved under edge flips  is slightly trickier to verify directly.  Namely, there are many cases to check, one such case is explicitly verified in figure \ref{boundary}.   By lemma \ref{baby23}, we can use a single triangulation per case to check all the possibilities.  However there are many cases.  To circumvent this issue we  can  boot strap from the  lower dimension examples.  
Namely, notice our needed result will follow  if the corresponding fact is true for the 2 dimensional supertetrahedra in figure \ref{suptets}.  Similarly, with this same reasoning, we find that the 2 dimensional fact is true provided the fact is true in one dimension.   In one dimension the result is transparent, see figure \ref{onedim}. 
This  inductive procedure can be easily suped up  to imply the analogs of these  results in all dimensions.

\begin{figure}
\center{\includegraphics{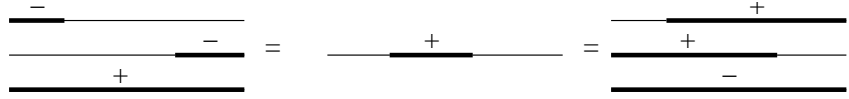}}
\caption{
\label{onedim}
On the left we see  a positively oriented one dimensional supersimplex and its 2 associated  one dimensional  prism halves, which are removed to form the simplex in the center of the figure. On the right, we see the consequence  of performing an edge flip, which results in the same  simplex, as needed.  }
 \end{figure}

\qed

\begin{figure}
\center{\includegraphics{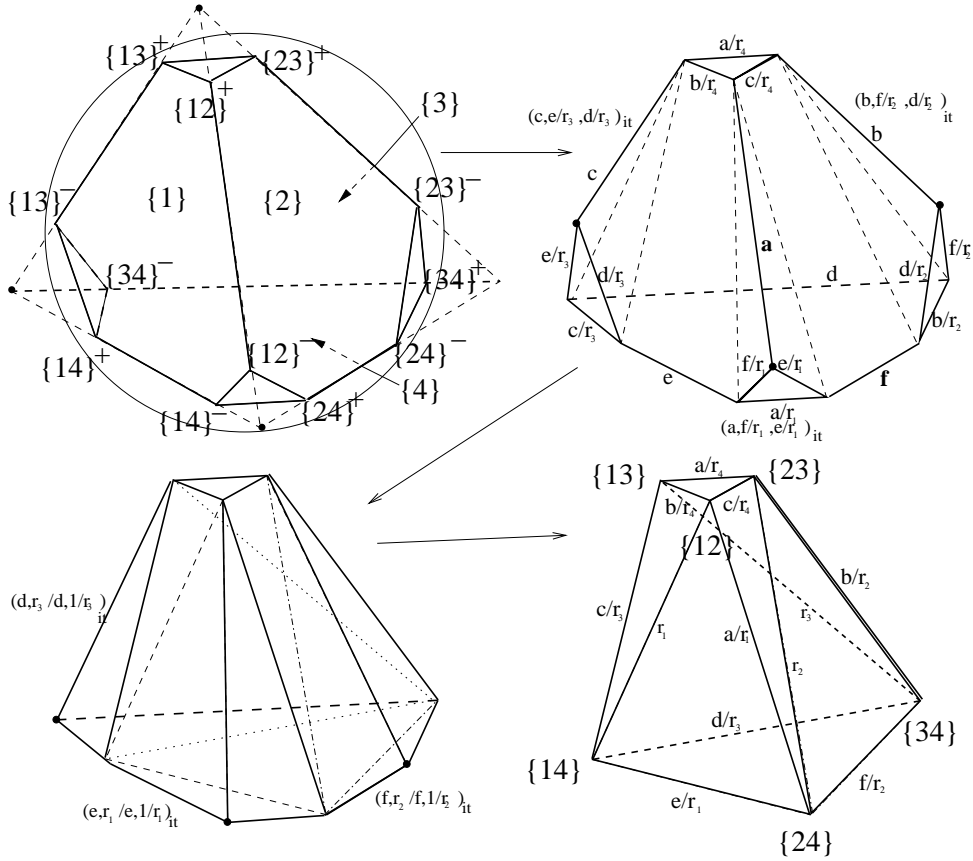}}
\caption{\label{heart}
Here we cut down to one of the supertetrahedron's $2^6$ hearts, namely the octahedron 
with the vertices
$ \{\{1,2\}^+,\{1,3\}^+\{1,4\}^+,\{2,3\}^+,\{2,4\}^+,\{3,4\}^+\}.$
This describes a $C$-region, $|H,R|$, determined from any supertetrahedron $|Su,R|$.
We have also labeled the clinants of  $|Su,R|$ and used these
 clinants to describe the clinants of $|H,R|$. 
 We have also indicated all the ideal tetrahedra utilized to cut   
 $|Su,R|$  down to $|H,R|$.    }
\end{figure}

{\bf A Supertetrahedron's Heart:} Given any convex ideal polyhedron,  $C$, if a vertex is $3$ valent,  then we can
{ cut off the ideal tetrahedron} it determines, hence,  determining a new convex ideal polyhedron in the process.  A {\it cut down} of $C$ , is a maximal sequence of such cuts.  For
example, cutting down  the convex prisms in figure \ref{prisms}   results in the empty set.
In figure \ref{heart},  we see an example of cutting down a 
convex supertetrahedron to an ideal octahedron. 
Once we have a cut down of a convex ideal polyhedron $C$, we may utilize our cutting down procedure 
 to  cut down any $C$-region $|C,R|$.
Cutting down a supertetrahedron always 
 result in some nonempty polyhedron,
 which we will call one of $C$'s {\it hearts}.   
 If we take the set of
 edges, $E=\{\{ij\}\}$, then  there are
$2^{|E|}$  hearts of a convex supertetrahedron each  given by the convex
hull of the  vertices 
 $\{ e^{sgn(e)}\}_{e \in E}$, where $sgn$ is any of the $2^{|E|}$ mappings in the form 
  \[ sgn:E \rightarrow \pm .\]

  In fact, this construction makes sense with respect to a supersimplex in every dimensions, where there are still $2^{|E|}$ such hearts, and any heart of any supersimplex  is
combinatorially equivalent to the polyhedron in Euclidean space
determined by taking the convex hull of the midpoints of the edges of the
$n$-simplex, called the abosimplex by Conway.  The ambotetrahedron happens to be the octahedron.

\begin{figure}
\center{\includegraphics{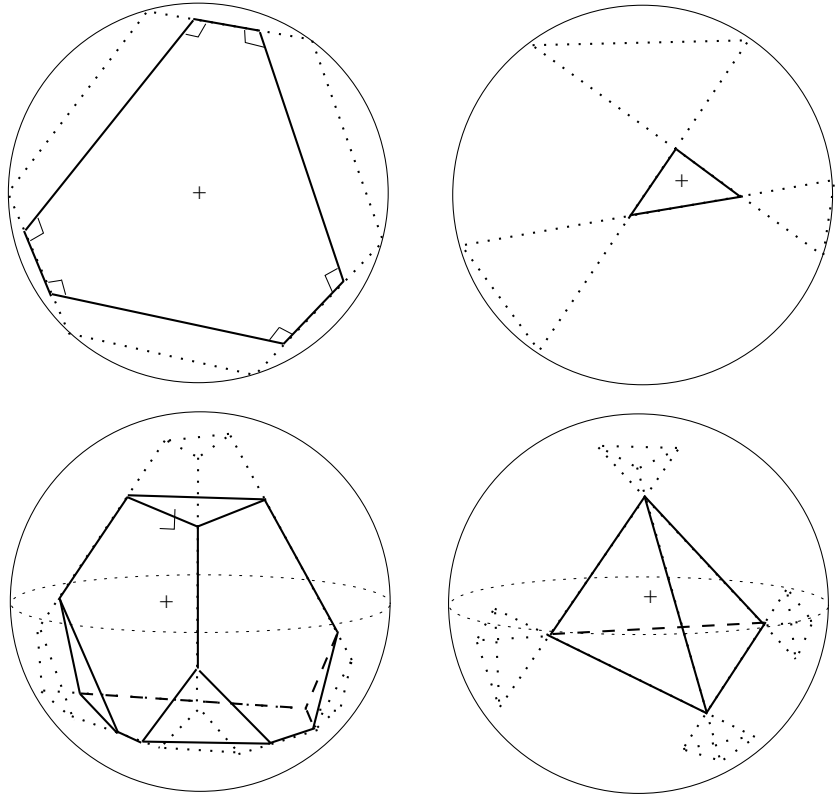}}
\caption{\label{tets}
 To define a generalized hyperbolic tetrahedron we use the fact from lemma \ref{prisl} that  each 
 $|P,R_{i}|$ has a well defined top half and we utilize 
the second part of lemma \ref{supl}, assuring us that we have a well defined region upon the removal from 
 $|Su,R|$ of the four top half prisms defined by the  $|P,R_{i}|$.  A 
 {\it  generalized hyperbolic tetrahedron} is the geometric region determined by 
removing the top half of each 
$|P,R_{i}|$. As in this figure, the notion of a 
generalized hyperbolic tetrahedron    coincides with the notion of a  finite
  hyperbolic tetrahedron, when each triple intersection of the tetrahedral planes  is non-empty. Let
$T(Su,R)$ denote the generalized  hyperbolic tetrahedron determined by
$|Su,R|$, and let $[T(Su,R)]$ be $T(Su,R)$'s  scissors class. }
\end{figure}

\begin{theorem}{\label{tetl}}
The following facts are true about  the 3 dimensional  generalized  hyperbolic tetrahedron as described in figure \ref{tets}.  
\begin{enumerate}
\item \label{p1}
A generalized  hyperbolic Tetrahedron  is  equivalent to a non-degenerate collection $4$
planes.
 \item \label{p3}
 \begin{eqnarray}\label{tetsis1}
2 [T(Su,R)] & = &  2[Su,R]-\sum_{i} [{P},R_i]. 
  \end{eqnarray}
   \begin{eqnarray}\label{tetsis2}
   2 [T(Su,R)]  & = & [Su,R]+[Su,\Omega(R)]. 
     \end{eqnarray}
    \begin{eqnarray}\label{tetsis3}
2 [T(Su,R)]  & = &   [H,R]+ [H,\Omega(R)].
     \end{eqnarray}
      \end{enumerate}
\end{theorem}
\proof
 Part \ref{p1} of this theorem follows from the second part of 
lemma \ref{supl} and the definition of the generalized  hyperbolic Tetrahedron.

From 
lemma \ref{prisl},  the top and
bottom halves of the  prism are scissors congruent.  
This together with the definition of the generalized  hyperbolic tetrahedron in figure \ref{tets}, gives us 
 formula \ref{tetsis1}. This fact is true in all dimensions. 

By lemma \ref{prisl},  $|P,R_i|$ is  isometric, via an orientation reversing
isometry, to $|P,\Omega(R)_i|$.  Hence, the half prisms used to define our generalized  hyperbolic tetrahedron 
cancel out as scissors classes, leaving us with  equation  (\ref{tetsis2}).   
Using the cut down described in figure \ref{heart}, we find the clinants of the ideal tetrahedron used to 
 cut   $|Su,R|$ down to  $|H,R|$ are conjugate  to those used to cut 
 $|Su,\Omega(R)|$ down to  $|H,\Omega(R)|$.
 Hence, by equation (\ref{neg}), the ideal tetrahedra utilized to cut down  $[Su,R]$ and $[Su,\Omega(R)]$ occur in oppositely oriented pairs, and cancel out as scissors classes. 
 Hence,   equation  (\ref{tetsis2}) implies  equation  (\ref{tetsis3}).  
 Formulas  (\ref{tetsis3}) and (\ref{tetsis2})  hold in all odd dimensions for the same reasons.
\qed 
 
{\bf Comment 1:}  Equation (\ref{tetsis1}) in theorem \ref{tetl} leads to interesting volume formula as well (see Mohanty \cite{Mo2} for a discussion of these volume formulas).
  Formula \ref{MY} will follow equation (\ref{tetsis3}), together with  the tools presented  in section \ref{group}.  Murakami and Yano present formula \ref{MY}  together with a pair of other formula (see \cite{Mu}) which follow from equation (\ref{tetsis1}),  together with  the tools presented  in section \ref{group}.
  
{\bf Comment 2:}  From theorem \ref{tetl} we nearly have theorem \ref{tetplane} from section \ref{com}.
 However, in theorem \ref{tetplane} the condition that each     $ \{ijk\} $  is ideal  has been dropped.  To extend the notion of generalized tetrahedra to include ideal points, one could simply drop the $R$'s injectivity condition, construct our supertetrahedra as just as in figure \ref{suptets}, and then  remove only the hyperideal and finite half prisms as in figure \ref{tets}.  We will arrive at  a well defined region which we can call a generalized hyperbolic   tetrahedron.  Such regions now include the possibility of  ideal vertices  and are in one to one correspondence with  collections of four, distinct, labeled, pairwise intersecting planes in  $H^3$, as need in theorem \ref{tetplane}.  
 There is a another version of this construction which is compatible with the constructions in sections \ref{clis} and \ref{group}.    To discuss this construction, one must go  outside the realm of $C$-regions altogether.  First,   observe that the union of a ideal tetrahedron  and its positively oriented mirror image  is a degenerate prism.  Such a union occurs  as the triangle in the waist of a convex  prism, see figure \ref{prisms},  degenerates to a Euclidean triangle. We can visualize this degenerations by fixing the positions on the sphere at infinity of the vertices at the bottom of the prism.  As we degenerate , we see a single ideal tetrahedron, the other becoming `hidden' at infinity. Such prisms are in the appropriate compactification of the space of prisms. Similarly, at every ideal vertex of a supertetrahedron there is a `hidden' ideal tetrahedra.   Removing our `hidden' half prisms will remove the supertetrahedron's hidden tetrahedra, hence, the resulting  tetrahedron  is still embedded and now has ideal points. 
 
  \end{subsection}

\begin{subsection}{\OP \label{clis}}
 
In this section, we will put coordinates on the space of
supertetrahedra, introduced in section \ref{simp}, and 
on the space of  ideal octahedra,  introduced in figure \ref{fireoct}. 
Technically,  when we  discuss coordinates on a space we will mean that an open dense set of the coordinates we present  form  the coordinates  of an open dense subset of our  space.   We shall also interpret mappings expressed in coordinates   as restricted to the appropriate open dense subsets.
We do this  because our $C$-regions are only designed to capture an open dense set of the objects of interest to us here, and because we will need to move between descriptions of our spaces that utilize complex coordinates and descriptions utilizing various clinants.    The need for such a convention can already be seen  when describing the space of ideal tetrahedra, see  the note at the end of section \ref{itet}.  In section \ref{dis}, we will carefully describe the exact open dense subsets to which the geometric constructions presented in this paper  apply without modification.

 \begin{figure}
\center{\includegraphics{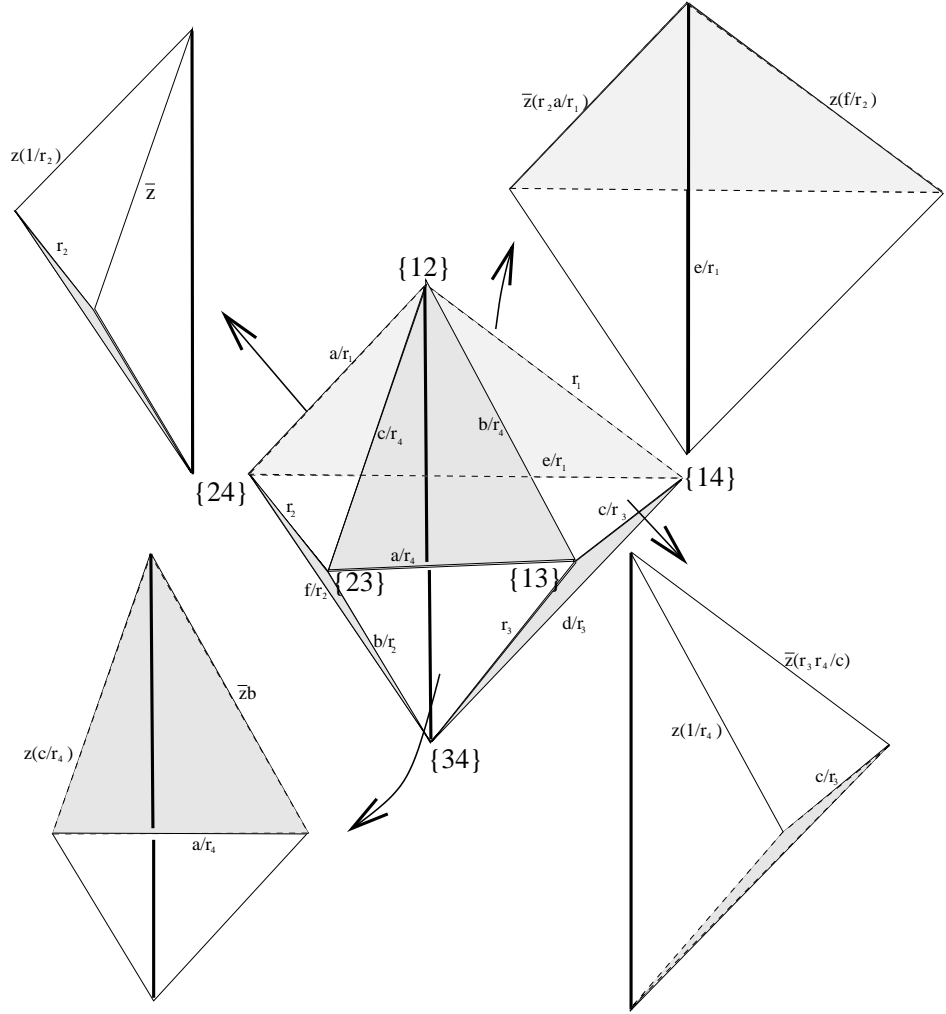}}
\caption{
\label{fireoct} Here we see a convex ideal octahedron, $O$,  in the
Klein model, with a specified triangulation, $S_o$.  Notice there are three such triangulations, each corresponding to choosing a pair of non-adjacent vertices. To each such triangulation there is also associated a {\it waist}, namely the four edges not containing either one of  the specified pair of non-adjacent vertices. 
  We let an ideal octahedron be  any $O$-region, $|O,R|$.
  If we wish to discuss the clinants of an ideal octahedron we 
let   $\theta^{ij}_{lk}$ denote the clinant associated to the edge with vertices  $\{ij\}$ and $\{lk\}$.
For example, in our figure   $\theta^{12}_{13} = \frac{b}{r_4}$.
  As pictured, the ideal  tetrahedral clinants of 
 $(S_o,R)$  are determined by  the octahedron's dihedral clinants
 up to a single unknown clinant, labeled  $z$. 
 }
\end{figure}

From figure \ref{fireoct}, we see that  an ideal octahedron can be  decomposed  into four ideal tetrahedra, hence 
  utilizing the complex coordinates of these ideal tetrahedra, we find that set of all
 \[ \vek{w} = (w_1 , w_4 , w_3  w_4  ) \in {\bf C}^4  \] 
  satisfying the  {\it holonomy constraint} that  
\begin{eqnarray}\label{hol}
 w_1w_2w_3w_4 & =  & 1
 \end{eqnarray}
form coordinates on the space of ideal octahedra. 
 We shall now attempt to utilize the octahedron's clinants to find another set of coordinates. 
To do so, note, the dihedral clinants of an ideal octahedron are $12$ unit complex number  indexed by the edges of the convex  octahedron that   multiply to $1$ at every vertex, and
that  multiply to $1$ around each of the octahedron's 3 waists. 
The indexing of the octahedral clinants is discussed in figure \ref{fireoct}, and 
 we shall let $\vek{o}$ denote an element of  $\times^{12} S^1$ that satisfies these octahedral constraints and is  ordered as follows 
  \[ \vek{o} = ( 
 \theta^{12}_{13},\theta^{12}_{14},\theta^{12}_{23},\theta^{12}_{24},
\theta^{34}_{13},\theta^{34}_{14},\theta^{34}_{23},\theta^{34}_{24},
\theta^{13}_{14},\theta^{13}_{23},\theta^{14}_{24},\theta^{23}_{24}
).   \] 
We shall see that such a $\vek{o}$ nearly determines an octahedron. 
To understand this claim it is useful to reparameterize our possible $\vek{o}$ via the 
coordinates introduced in figure \ref{heart}.
To accomplish this, first notice from figure \ref{heart}, we see that  for a suitable  $(a,b,c,d,e,f;r_1,r_2,r_3,r_4)$ that 
that we can form an octahedron  with clinants 
 $  \vek{o}(a,b,c,d,e,f;r_1,r_2,r_3,r_4)$ equal to 
  \[   
\left(\frac{b}{r_4}, r_1, \frac{c}{r_4},\frac{a}{r_1},r_3,\frac{d}{r_3},\frac{b}{r_2},\frac{f}{r_2},\frac{c}{r_3},\frac{a}{r_4},\frac{e}{r_1},r_2 \right) . 
\]
 Let us comment on these  $(a,b,c,d,e,f;r_1,r_2,r_3,r_4)$ coordinates. 
 To each $\{ij\}$ edge  in figure \ref{heart} we have associated a clinant. 
 The fact that the clinants at an ideal vertex multiply to one determine the $r_i$ up to sign. Namely,
 \begin{eqnarray}
 r_1^2&  = &aef  \label{s1} \\
 r_2^2 & =& bdf   \label{s2} \\ 
 r_3^2 & = & c de   \label{s3} \\
 r_4^2 & =& abc.  \label{s4} 
 \end{eqnarray}
 There is also a condition on the $r_i$ due to  the fact that the clinants around  an octahedron's waist multiply to one; namely \begin{eqnarray}\label{s5}
r_1 r_2 r_3 r_4 &  = &  abcdef.
\end{eqnarray}
We will let  a coordinate in the form 
\[ \vek{s} = (a,b,c,d,e,f;r_1,r_2,r_3,r_4)   \] 
denote an element of   $\times^{10} S^1$ that 
 satisfies conditions         \ref{s1}-\ref{s5}. 
These can be used as a reparameterization of the octahedron's  clinants since the mapping 
\begin{eqnarray*}\label{supc}
\vek{s}(\vek{o}) & =  & \left(
 \theta^{12}_{24} \theta^{12}_{14},
  \theta^{23}_{34} \theta^{23}_{24},
\theta^{13}_{14} \theta^{13}_{34},
\theta^{14}_{34} \theta^{13}_{34},
\theta^{14}_{24} \theta^{12}_{14},
\theta^{24}_{34} \theta^{23}_{24},
\theta^{12}_{14},\theta^{23}_{24},\theta^{13}_{34},\frac{1}{\theta^{12}_{13}\theta^{13}_{23}\theta^{12}_{23}}
 \right).
 \end{eqnarray*}
is easily checked to be  $\vek{o}$'s inverse.

In figure \ref{fireoct}, we see that the space of octahedra can be constructed from 
 dihedral clinants  once we have determined  the  clinant labeled  $z$.
If we let
  \begin{eqnarray}
 (\vek{s},z)^{\mat{}}_{\sr{M}\sr{Y}} & =  & \left(
\frac{1}{z},
\frac{z}{ r_2},
\frac{b}{z},
\frac{zc}{r_4},
\frac{ r_3 r_4}{zc},
\frac{z}{ r_4},
\frac{  a r_2}{ z r_1}
\frac{z f}{ r_2}
 \right)
  \end{eqnarray}
  and let $(\vek{s},z)^{\mat{}}_{\sr{M}\sr{Y}}(j)$ be the  $j^{th}$ component of $(\vek{s},z)^{\mat{}}_{\sr{M}\sr{Y}}$, then, from equation (\ref{fund}), we have that the  $w_i$ coordinates of our  ideal octahedron are equal to 
  \[{w}_{i}(\vek{s},z) = \frac{1- \overline{(\vek{s},z)^{\mat{}}_{\sr{M}\sr{Y}}(2i)} }{1-  (\vek{s},z)^{\mat{}}_{\sr{M}\sr{Y}}(2i-1)   }. \]
Notice if we let 
\begin{eqnarray}
 (\vek{s})_{waist} & =  & \left(
r_2,
\frac{a}{r_4},
\frac{c}{r_3},
\frac{e}{r_1}
\right)
  \end{eqnarray}
then by equation (\ref{mis}), for  any unit sized $z$, 
 \begin{eqnarray}\label{waist}
 2 [\vek{w}(\vek{s},z)] & = & [ (\vek{s},z)^{\mat{}}_{\sr{M}\sr{Y}}] + [ (\vek{s})_{waist}] 
 \end{eqnarray}
 
To find the $z$ that corresponds to our octahedron, notice that 
 equation (\ref{hol}),  implies that  $z$ must solve 
 \begin{eqnarray}\label{hol2}
k(\vek{s},z) & = & \Pi_{i=1}^4(1-  \overline{(\vek{s},z)^{\mat{}}_{\sr{M}\sr{Y}}(2i)} ) -   \Pi_{i=1}^4(1-  (\vek{s},z)^{\mat{}}_{\sr{M}\sr{Y}}(2i-1)  ) =0.
\end{eqnarray}  
Upon  multiplying  out this equation (\ref{hol}), 
we find that $(1/z)k(\vek{s},z)$ is a quadratic polynomial in $z$ with coefficients in the  supertetrahedral 
clinants. Hence an octahedron 
is determined by  its $\vek{s}$ coordinate  along with the correct root of $k(\vek{s},z)$.  
We will call the root, $\mat{\rho}$, corresponding to our octahedron the octahedron's {\it octahedral root}.  
In order, to understand the space of octahedron of fundamental importance is  the following
quantity: 
\[ \mat{\delta} = \frac{-4}{ab^2d} disc\left( \frac{k(\vek{s},z)}{z} \right), \] 
where $disc$ refers to the quadratic's discriminant.

\begin{lemma}\label{octlem}
For any $\vek{s}$ we have that $\mat{\delta}$ is real and the  open set of $(\vek{s},\mat{\rho})$ with $\mat{\delta}>0$ form coordinates on the space of octahedra.
\end{lemma}
\proof
Let   $\vek{c} = (\sr{A},\sr{B},\sr{C},\sr{D},\sr{E},\sr{F}) \in \times^6 S^1$ and define  
\[  \vek{s}(\vek{c}) = (\sr{A}^2, \sr{B}^2,\sr{C}^2,\sr{D}^2,\sr{E}^2,\sr{F}^2,-\sr{A}\sr{E}\sr{F},-\sr{B}\sr{D}\sr{F},-\sr{C}\sr{D}\sr{E},-\sr{A}\sr{B}\sr{C}  ).\]
Via this mapping  the $\vek{c}$ coordinates form a $2^3$ fold cover of our    $\vek{s}$ coordinates. 
Let $\negg$ be the set of  transformations that  negate a collection of the  
$\vek{c}$ coordinates that share a face in figure   \ref{label}. 
This cover's  deck group is the copy of $\frac{\Bbb{Z}}{8 \Bbb{Z}}$  generated by the  elements of $N$.
As a set, this deck group is $\trivial \bigcup \negg$, where $\trivial$ is the group introduced in the comment at the end of section \ref{com}.
We find that $k(\vek{s}(\vek{c}),z)$
 can be expressed as 
  \[\mat{h}(\vek{c},z)= \frac{2}{ \sr{A} \sr{B}^2 \sr{D} z} k(\vek{s}(\vek{c}),z) ={\mat{\alpha}} z^2 + 2 \mat{\beta} z + \bar{\mat{\alpha}} \]
  with
  \[ \bar{\mat{\alpha}}(\vek{c}) =
2 \left( \frac{\sr{B}}{\sr{E}}+\sr{A}\sr{D}+\sr{B}\sr{E}+\sr{A} \sr{B}^2 \sr{D}  + 
\sr{A}\sr{B}\sr{F}+\frac{\sr{A}\sr{B}}{\sr{F}}+\sr{B}\sr{C}\sr{D}+\frac{\sr{B}\sr{D}}{\sr{C}}   \right) \]
and 
\[ \mat{\beta}(\vek{c})  = \frac{\sr{A}}{\sr{D}}+\frac{\sr{D}}{\sr{A}}+\frac{\sr{C}}{\sr{F}}+\frac{\sr{F}}{\sr{C}}+\sr{C}\sr{F}+\frac{1}{\sr{F}\sr{C}}-\left(\sr{A}\sr{D}+\frac{1}{\sr{A}\sr{D}}+\frac{\sr{B}}{\sr{E}}  +\frac{\sr{E}}{\sr{B}} +\sr{B}\sr{E}  +\frac{1}{\sr{B}\sr{E}} \right),
 \] 
and notice $\beta$ is real.
These quantities are not independent of the choice of $\vek{c}$ satisfying $\vek{s}(\vek{c})= \vek{s}$.
 Namely, $\mat{\alpha}(D \cdot \vek{c}) = -\mat{\alpha}(\vek{c})$, and that $\mat{\beta}(D \cdot \vek{c}) = -\mat{\beta}(\vek{c})$ when $D \in N$,
while the other  deck transformations preserve these quantities. 
 In particular, 
 \[ \mat{\delta}(\vek{c})  = |\mat{\alpha}|^2 - \mat{\beta}^2, \]
is real and dependent only on $\vek{s}$.
Furthermore when  
 $\mat{\delta} >0$ we have that 
 \[ \mat{\rho}(\vek{c})  =  \frac{-\mat{\beta}  - i  \sqrt{\mat{\delta}} }{{\mat{\alpha}}}. \]
 and $\mat{\rho}(D \cdot \vek{c})$  are the two root of $\mat{h}(\vek{c},z)$.
 Both these roots are 
 unit sized since 
\[ |\mat{\rho}|^2 =  \frac{\mat{\beta}^2+\mat{\delta}}{\mat{\alpha} \bar{\mat{\alpha}}} = \frac{|\mat{\alpha}|^2}{|\mat{\alpha}|^2}=1. \] 
 Hence when $\mat{\delta} >0$  we can generically construct our needed ideal tetrahedra.
Now when $\mat{\delta} <0$
we have that $\mat{h}(\vek{c},z)$'s roots are 
given by 
\[ \frac{-\mat{\beta}  \pm    \sqrt{\mat{\delta}} }{{\mat{\alpha}}}, \] 
 and, hence,  have magnitude 
\[  \frac{\mat{\beta}^2+\mat{\delta}-2\beta \sqrt{\mat{\delta}}}{\mat{\alpha} \bar{\mat{\alpha}}}   = 1 \pm \frac{2\beta \sqrt{\mat{\delta}}}{|\mat{\alpha}|^2}, \]
which is not  not unit sized, since
  $\delta <0$ implies that   $\beta^2 >0$, hence for 
$ \beta \sqrt{\mat{\delta}}$ to be equal to $0$, we need that $\delta=0$.

\qed

We will freely replace the coordinate $(\vek{s},\mat{\rho})$ of an ideal tetrahedron 
with a $(\vek{c})$,  as introduced in the proof of lemma \ref{octlem},
where $\vek{s}(\vek{c})=\vek{s}$ and
$\mat{\rho}(\vek{c})=\mat{\rho}$.
We will say that $(\vek{c})$  is a choice of  {\it tetrahedral circulants} associated to $(\vek{s},\mat{\rho})$. 
Notice the tetrahedral circulants form a $4$ fold cover of the  $(\vek{s},\mat{\rho})$ with nontrivial deck transformations 
given by the transformations that negate all the tetrahedral clinants except those corresponding to  an opposite pair of edges. 
These are the geometrically trivial transformations which arose in the sketch of theorem \ref{res2}'s proof. 
Let  $(\vek{c})_{oct}$ denote the octahedron with tetrahedral circulants  $(\vek{c})$.
Let 
\[ (\vek{c})_{op}= \left((\vek{c})_{oct},\left(\bar{\vek{c}}\right)_{oct}\right) \]
be called an pair of  {\it \op}.

\begin{lemma}\label{sign}
The $(\vek{s},\mat{\rho})$ with $\mat{\delta}>0$ form coordinates on the space of all \op,
\begin{eqnarray}\label{opii}
2  [(\vek{c})_{op} ]
 & = & 
  [ (\vek{c})^{\mat{}}_{\sr{M}\sr{Y}}] + \left[ \left(\bar{\vek{c}}\right)^{\mat{}}_{\sr{M}\sr{Y}}\right],  
\end{eqnarray}
\begin{eqnarray}\label{sn3}
 \left[\left(\bar{\vek{c}}\right)_{op}\right] & = &  [ (\vek{c})_{op}], 
 \end{eqnarray}
and if $D \in \negg$, then  $D$ commutes with conjugation, 
\begin{eqnarray}\label{sn2}
 [ (D \cdot \bar{\vek{c}})_{oct}]    & = &  - [(\vek{c})_{oct}], 
 \end{eqnarray}
 and 
 \begin{eqnarray}\label{sn4}
 [ (D \cdot \vek{c})_{op}]    & = &  - [(\vek{c})_{op}]. 
 \end{eqnarray}
\end{lemma}
\proof 
That the  $(\vek{s},\mat{\rho})$ with $\mat{\delta}>0$ form coordinates on the space of \op follows immediately from lemma \ref{octlem}.  Equation (\ref{opii}) follows from equation (\ref{waist}) and equation (\ref{neg}). 
Equation (\ref{sn3}) follows from the  definition of $(\vek{c})_{op}$.  We mention it  in order   to  emphasize that conjugation does {\bf not} correspond to reversing orientation.
Instead $D \in \negg$ does the job of reversing orientation job, even on the level of the octahedron. 
To see this notice that if we conjugate every clinant in figure \ref{fireoct} then we arrive at  
 $(D \cdot \bar{\vek{c}})_{oct}$, since 
 $ \mat{\rho}(D \cdot \bar{\vek{c}})=\overline{{\mat{\rho}}\left( \bar{\vek{c}}\right)} $.
 Hence, equation (\ref{sn2}) follows form this observation together with equation (\ref{neg}).
 The equation  (\ref{sn4}) follows from  equations (\ref{sn3}) and (\ref{sn2}).
 
 \qed

Generically,  we can  geometrically invert an octahedron back into a supertetrahedron utilizing the ideal tetrahedra in figure \ref{heart}.    
 In particular, by lemma \ref{octlem},  the coordinates  $(\vek{c})$ with $\mat{\delta}>0$  will cover the space of  generalized hyperbolic tetrahedron.  We will denote the tetrahedron corresponding to  such a  coordinate as   $(\vek{c})_{tet}$.

\begin{lemma}\label{main}
Generically,
\begin{eqnarray}\label{siso}
 2[(\vek{c})_{tet}] &  = &  [(\vek{c})_{op}].
\end{eqnarray}
\end{lemma}
\proof
Let $|Su,R|$ be the supertetrahedron generically  corresponding to  $(\vek{c})$, and so
\[  (\vek{c})_{oct}  = |H,R|. \]
 Since $|H,\Omega(R)|$'s is an octahedron with  dihedral clinants  $\vek{ \overline{d}(\vek{c})}$,
 either $ |H,\Omega(R)|  =  \left(\bar{\vek{c}}\right)_{oct}$ or  $ |H,\Omega(R)|  =  \left(D \cdot \bar{\vek{c}}\right)_{oct}$ 
 with $D \in \negg$.
If    $ |H,\Omega(R)|  =  \left(D \cdot \bar{\vek{c}}\right)_{oct}$ 
 Then $ |H,\Omega(R)|$ has all its tetrahedral clinants conjugate to those of $|H,R|$, and hence by equation (\ref{neg}),   
\[ [H,R]+[H,\Omega(R)] = 0.\]
This, generically, contradicts equation (\ref{tetsis3}), 
hence  $ |H,\Omega(R)|  =  \left(\bar{\vek{c}}\right)_{oct}$ and 
    \[ (|H,R|, |H,\Omega(R)|) = (\vek{c})_{op}. \]

\qed

  {\bf Comment:}  Notice that when $\delta = 0$, that we have a unique unit root
  $\-\mat{\beta}/
  {\mat{\alpha}}$ of 
  $\mat{h}(\vek{c},z)$.  Hence are still in a position to construct our octahedra.  In this case, we find that all our 4 tetrahedral planes intersect in a point, which  corresponds to an infinitesimal hyperbolic tetrahedron, or rather a Euclidean tetrahedron.  One nice way to understand this is to note that $\delta = -16 \det(Gr)$, where $Gr$ is the Graham matrix associated to our planes.  Hence  $\delta=0$ exactly when we are in the Euclidean case.

 \end{subsection}

\begin{subsection}{The  Group  \label{group}}

\begin{figure}
\center{\includegraphics{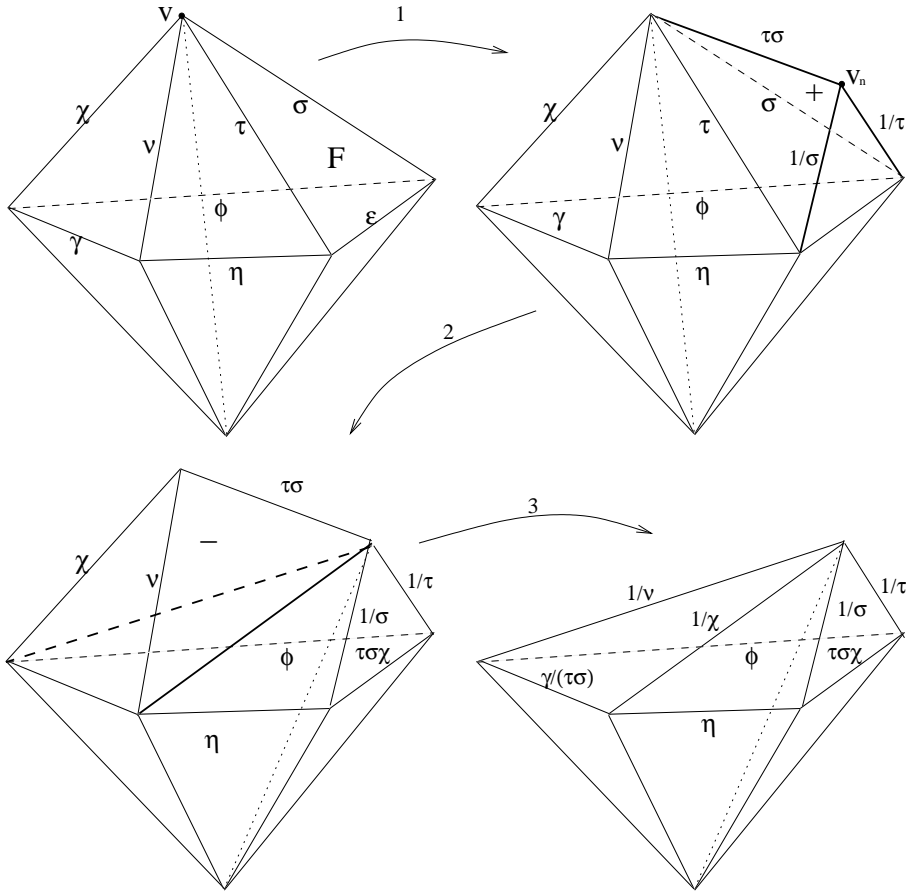}}
\caption{
\label{puffncut}
We start by describing the   puff of   $|O,R|$ with
respect to a vertex,  $v$, and face, $F $,  with $v \in F$. 
If we assume  
$\tau \sigma \neq 1$,  then $F$'s neighboring faces at $v$
intersect in a geodesic, pictured as the edge  with the   $\tau \sigma$ clinant attached to
it after step 1.  Using this geodesic we may determine a new ideal
vertex, the vertex labeled
$v_n$.   A {\it puff} is the operation of adjoining   the ideal tetrahedron
$(1/\tau,1/\sigma,\tau \sigma)_t$  that contains $v_n$ and $F$, as pictured after step 1.
After our puff, $v$ becomes a trivalent vertex, as  in step 2. 
Being a  trivalent vertex,   we may {\it  cut} $v$ off, as seen in step 3.
 Equivalently, we may change $R$ to $P^{F}_v R$ by demanding that
$P^{F}_v R(p)=v_n$ for the for the $p$ such that 
 $R(p)=v$, and letting $R(q)=P^{F}_v R(q)$ for every other octahedral vertex  $q$. 
 We   say   $|O,P^{F}_v R|$  results from $|O, R|$ via a {\it puff-and-cut}.
  }
\end{figure}

   We introduce the group described in theorem \ref{res2}.   It will be generated by the puff-and-cuts in 
in figure \ref{puffncut}.   These puff-and-cuts  take an octahedron  $|O, R|$ and transform it into 
 $|O,P^{F}_v R|$. For every vertex face pair this induces  a mapping of the octahedral clinants, 
which we shall dente as $P_{v}^F \vek{o}$.  
For example, letting $F$ be the face of our octahedron with $F=\{\{12\},\{13\},\{23\}\}$
 we have  $P_{\{12\}}^F \vek{o}$ equals 
\[   P_{\{12\}}^F \vek{o} = \left( 
\frac{1}{\theta^{12}_{14}},\frac{1}{\theta^{12}_{13}},\frac{1}{\theta^{12}_{24}},\frac{1}{\theta^{12}_{23}},
\theta^{34}_{13},\theta^{34}_{14},\theta^{34}_{23},\theta^{34}_{24},
\theta^{13}_{14} \theta^{12}_{13} \theta^{12}_{14},\theta^{13}_{23},\theta^{14}_{24}
,\frac{ \theta^{23}_{24}}{\theta^{12}_{13} \theta^{12}_{14}}
\right).   \] 
We can express this transformation in the supertetrahedral coordinates as 
\[P_{\{12\}}^F \cdot  \vek{s} 
= \vek{s}(P_{\{1,2\}}^F (\vek{o}(\vek{s}))) 
=  \left( \frac{1}{a},b,c,d,e,f,\frac{r_1}{a},r_2,r_3,\frac{r_4}{a} \right). \]
Let $\wacky$ denote the group generated by the $P_{v}^F$, as in the sketch of theorem \ref{res2}'s proof.
Geometrically, we will be most interested in the action of $\wacky$ when viewed as 
acting on \op.    We explore  $\wacky$'s  algebraic structure  in section \ref{algsym}, for now we need the following lemma.

\begin{lemma}\label{octplem}
For $\gen \in G$ we generically  have 
 \begin{eqnarray}\label{octplem2}
 [(\vek{c})_{op}] &=&    [(\gen \cdot \vek{c})_{op}] ,
  \end{eqnarray}
   \begin{eqnarray}\label{commute}
 \left(\overline{\gen \cdot \vek{c}}\right)   & = &
 \left({\gen \cdot \bar{\vek{c}}}\right).
  \end{eqnarray}
 \end{lemma}
 \proof
Notice that  the clinants of 
 \[  (\vek{c})_{oct}=|H,R| \] 
 and 
 \[  \left(\bar{\vek{c}}\right)_{oct}=|H,\Omega(R)| \]
  are conjugates, hence the ideal tetrahedra utilized in figure \ref{puffncut} to transform 
   $ (\vek{c})_{op}$
   to 
   \[  (|H, P_{v}^F R|, |H, P_{v}^F \Omega(R)|) = (P_{v}^F \cdot \vek{c})_{op}  \]
   will  cancel out as scissors classes,  by equation (\ref{neg}).

To prove equation (\ref{commute}), note       we either have equation \ref{commute} or
  \[
 \left({\gen \cdot \bar{\vek{c}}}\right) =  \left(D \cdot \overline{\gen \cdot \vek{c}}\right)    \] 
 for $D \in \negg$.  Notice  by equation  \ref{sn4} 
 that 
  \[ [\left(D \cdot \overline{\gen \cdot \vek{c}}\right)_{op} ] = - [\left(\overline{\gen \cdot \vek{c}}\right)_{op}] . \]
  By the definition of $\gen \in G$ and equation \ref{sn3} we have 
    \begin{eqnarray*}
   [\left({\gen \cdot \bar{\vek{c}}}\right)_{op}] & = & [\gen \cdot \left(\bar{\vek{c}}\right)_{op}]  \\
   & = & [\gen \cdot \left(\vek{c}\right)_{op}]   \\
  & = &  [\left( \gen \cdot \vek{c}\right)_{op}]  \\
   & = &   [\left( \overline{\gen \cdot \vek{c}} \right)_{op}] 
  \end{eqnarray*}
 Hence, $\left( \overline{\gen \cdot \vek{c}} \right) \neq \left(D \cdot \overline{\gen \cdot \vek{c}}\right)$ and    
 equation (\ref{commute}) follows.

\qed

The 24 possible puff-and-cuts in figure \ref{puffncut} 
 are the generators of $\wacky$.  It is straight-forward to verify that the opposite
faces at a vertex induces the same transformation of the $\vek{s}$ coordinate.
From this, these  transformations    induce the same transformation of the  \op, due the following 
 scholium to lemma \ref{octplem}.

\begin{lemma}\label{prep}
How  $\gen \in G$ acts  on $(\vek{c})_{op}$
is completely determined by how $\gen$ acts on 
$\vec{d}(\vec{t})$.   
\end{lemma}
\proof
The $\vek{c}$ coordinate  representing 
 $\gen  \cdot \vek{s}$. 
 is either 
$\gen \cdot \vek{c}$ or $D \cdot \gen \cdot \vek{c}$ with $D \in \negg$.
 But 
  \[ [(\gen \cdot \vek{c})_{op}] =- [(D \cdot \gen \cdot \vek{c})_{op}]\] 
 hence, by lemma \ref{octplem}, generically only one of $\gen \cdot \vek{c}$ or $D \cdot \gen \cdot \vek{c}$, 
  can represent a   image of  $(\vek{c})_{op}$  under $\wacky$.

\qed

  Hence, $\wacky$ is generated by twelve 
 of the $P^F_v$.   These generators break naturally into two class,
those with a {\it shaded face}  in
figure \ref{fireoct}, and those with an {\it unshaded} face in
figure \ref{fireoct}.  Hence we will simply denote $P_{v}^F $  as $P_{v}^u$ or  $P_{v}^s$  depending on whether we are using one of the shaded or unshaded faces at $v$. 
For example, our above $P_{\{1,2\}}^F =  P_{\{1,2\}}^s$, since we utilized a shaded face.
 As with our $P_{\{1,2\}}^s$ transformation, in the tetrahedral coordinates
  all of our shades transformation simply invert the clinant corresponding  to $v$. 
In terms of our supertetrahedral coordinates, we have 
\begin{eqnarray}\label{g1}
P_{\{12\}}^F \cdot  \vek{s}  
=  \left( \frac{1}{a},b,c,d,e,f,\frac{r_1}{a},r_2,r_3,\frac{r_4}{a} \right)
\end{eqnarray}
\begin{eqnarray}\label{g2}
 P_{\{34\}}^s \cdot\vek{s}
=  \left(a,b,c,\frac{1}{d},e,f,r_1,\frac{r_2}{d},\frac{r_3}{d},r_4 \right)
\end{eqnarray}
\begin{eqnarray}\label{g3}
P_{\{23\}}^s \cdot\vek{s}
=  \left(a,\frac{1}{b},c,d,e,f,r_1,\frac{r_2}{b},r_3,\frac{r_4}{b} \right) 
\end{eqnarray}
\begin{eqnarray}\label{g4}
 P_{\{14\}}^s \cdot\vek{s}
=  \left(a,b,c,d,\frac{1}{e},f,\frac{r_1}{e}r_2,\frac{r_3}{e},r_4 \right)
\end{eqnarray}
\begin{eqnarray}\label{g5}
P_{\{13\}}^s \cdot\vek{s}
=  \left(a,b,\frac{1}{c},d,e,f,r_1,r_2,\frac{r_3}{c},\frac{r_4}{c} \right) 
\end{eqnarray}
\begin{eqnarray}\label{g6} P_{\{2,4\}}^s \cdot\vek{s}
=  \left(a,b,c,d,e,\frac{1}{f},\frac{r_1}{f},\frac{r_2}{f},r_3,r_4 \right). 
\end{eqnarray}
These shaded elements  generate the group isomorphic to  $\left( \frac{\bf Z}{2 {\bf Z}} \right)^6$ 
discussed in proof of corollary \ref{cor}, which we will call the {\it shaded subgroup}.
If our octahedron is $|H,R|$ for some supertetrahedron $|Su,R|$, then   $P_{v}^s$ applied to our octahedron corresponds to $|H,F_{v} R|$.  In other words,  $P_{v}^s$  corresponds to flipping the $v$ edge of $|Su,R|$.
In particular, by theorem \ref{tetl}, the shaded subgroup  preserves not only the scissor class of the tetrahedron but the underlying tetrahedron itself.  In fact, recall from section \ref{clis} that the  space of \op covers  the generalized hyperbolic tetrahedra.  The shaded subgroup is the group of deck transformations of this cover. 

 The unshaded transformations are much more geometrically subtle. Algebraically they are  given by the following transformations.  
\begin{eqnarray}\label{g7}
  P_{\{12\}}^u \cdot\vek{s}
=  \left(a,\frac{r_4}{r_1},\frac{r_1r_4}{a},d,\frac{er_4}{br_1},\frac{fr_4}{br_1},\frac{r_4}{b},\frac{r_3}{e},r_3,r_4  \right)
\end{eqnarray}
\begin{eqnarray}\label{g8} 
P_{\{34\}}^u \cdot\vek{s}
=  \left(a,\frac{r_2}{r_3},  \frac{c r_2}{b r_3},d,\frac{e r_2}{b r_3},\frac{r_2 r_3}{d}, r_1,r_2, \frac{r_2}{b}, 
 \frac{ r_1}{e}\right) 
 \end{eqnarray}
\begin{eqnarray}\label{g9}
P_{\{23\}}^u \cdot\vek{s}
 = \left(  \frac{r_2r_3}{b}, b, \frac{r_4}{r_2}, \frac{d r_4}{c r_2},e,\frac{fr_4}{cr_2},r_1, \frac{r_4}{c},\frac{r_1}{f}, r_4 \right) 
 \end{eqnarray}
\begin{eqnarray}\label{g10}
 P_{\{14\}}^u \cdot\vek{s}=  \left(\frac{ar_3}{cr_1}, b, \frac{r_3}{r_1}.\frac{r_1r_3}{e},e,\frac{fr_3}{cr_1}.\frac{r_3}{c}.r_2,r_3,\frac{r_2}{ f} \right) 
 \end{eqnarray}
\begin{eqnarray}\label{g11}
P_{\{13\}}^u \cdot\vek{s}
   = \left( \frac{r_4}{r_3},  \frac{r_3 r_4}{c}, c,\frac{d r_4}{a r_3},\frac{e r_4}{a r_3},f,\frac{r_2}{d} ,r_2,\frac{r_4}{a},r_4 \right) \end{eqnarray}
\begin{eqnarray}\label{g12}
 P_{\{24\}}^u \cdot\vek{s}
=  \left(\frac{r_1}{r_2},  \frac{b r_1}{a r_2}c,\frac{d r_1}{a r_2},\frac{r_1 r_2}{f},f, r_1, \frac{r_1}{a},r_3,\frac{r_3}{ d} \right) 
\end{eqnarray}

 \begin{subsubsection}{Proof of Theorem \ref{fundsis} \label{funpr}}
 
 We can now use $\wacky$ to prove theorem 
\ref{fundsis}.  
Namely, we use the fact that $\osa$ from the sketch of the proof of theorem \ref{fundsis} is in $\wacky$. 
In fact,  if we  let $s_1$ and $s_2$ be the  elements of the shaded subgroup determined by 
 \[ s_1 \cdot\vek{s} =  \left(\frac{1}{a},\frac{1}{b},c,\frac{1}{d},\frac{1}{e},f,
 \frac{r_1}{ae},\frac{r_2}{bd},\frac{r_3}{de},\frac{r_4}{ab} \right) \]
 \[ s_2 \cdot\vek{s} =  \left(a,b,\frac{1}{c},d,\frac{1}{e},f,
 \frac{r_1}{e},r_2,\frac{r_3}{ce},\frac{r_4}{c} \right). \]
then
\[\osa = s_2   P_{\{3,4\}}^u s_1. \]
Note 
 \[ \osa \cdot\vek{s} =
 \left( \frac{1}{a}, \frac{d}{r_2r_3},\frac{e r_2}{b r_3}, \frac{1}{d}, \frac{b r_3}{c r_2}, \frac{r_2}{r_3},
 \frac{b}{r_4}, \frac{1}{r_3}, \frac{e}{r_3},\frac{1}{r_4} \right). \]
 
  We chose $\osa$ so that we could take the hat off  $\mat{\rho}$ from section \ref{clis}. 
  In other words,  we have the following lemma.

  \begin{lemma}\label{inter2}
 \begin{eqnarray*}
   \rho(\vek{c}) & = &  \mat{\rho}(\osa \cdot \vek{c})  
\end{eqnarray*}
\end{lemma} 
 \proof
This nearly follows by direct substitution.
The nearly refers to the fact that 
by direction substitution we find for each $\vek{c}$ that 
 $    \rho(\vek{c}) = \mat{\rho}(\osa \cdot \vek{c})  $
or 
 $    \rho(\vek{c}) = \mat{\rho}(\osa \cdot (D \cdot  \vek{c}))  $, 
 with $D \in \negg$.
 By continuity of the construction with respect to the parameters, 
 we need to verify the formula 
 for an element of each of the  connected components of the subset of
 $\vek{c}$ coordinates where $\mat{\delta} >0$.  
 There are only two such components which are related by the transformation sending 
 $\vek{c}$ to  $D \cdot  \vek{c}$, and the formula is easily verified.
 
 \qed

 From lemma \ref{inter2} and a direct substitution we have that 
  \begin{eqnarray}\label{inter}
 (\osa  \cdot \vek{c})^{\mat{}}_{\sr{M}\sr{Y}}&  =  & (\vek{c})_{\sr{M}\sr{Y}}. 
\end{eqnarray}

  We may now prove theorem   \ref{fundsis},
by noting that  generically
\[ \begin{array}{lllll}
 4 [(\vek{c})_{tet}] &  =  &2  [(\vek{c})_{op}]  & \mbox{by} & \mbox{equation } \ref{siso}  \\
& = &  2 [(\osa \cdot \vek{c} )_{op}]   & \mbox{by} & \mbox{equation }  \ref{octplem2}   \\ 
& = & [ (\osa  \cdot \vek{c})^{\mat{}}_{\sr{M}\sr{Y}}] + 
\left[ \left(\overline{\osa \cdot \vek{c}}\right)^{\mat{}}_{\sr{M}\sr{Y}} \right]  & \mbox{by} & \mbox{equation }  \ref{opii}\\
& = &  [ (\osa  \cdot \vek{c})^{\mat{}}_{\sr{M}\sr{Y}}] +  \left[ \left({\osa \cdot \bar{\vek{c}}}\right)^{\mat{}}_{\sr{M}\sr{Y}} \right]  & \mbox{by} & \mbox{equation }  \ref{commute}\\
&=& [(\vek{c})_ {\sr{M}\sr{Y}} ] + \left[\left(\bar{\vek{c}}\right)_{\sr{M}\sr{Y}}  \right] & \mbox{by} & \mbox{equation }  \ref{inter}.
\end{array} \]

 \end{subsubsection}

 \begin{subsubsection}{An Algebraic Description of $\wacky$ \label{algsym}}

     To  expose the algebraic structure of $\wacky$  it is useful to  introduce some new coordinates. 
Recall that to resolve the branching necessary to determine $\mat{\rho}$ from  $\vek{s}$ we only need a two fold cover of $\vek{s}$, not the $2^3$ fold cover determined by the $\vek{c}$ coordinates.  We now describe this cover. 
  Let $\vek{b} = (t,u,v,T,U,V;r) \in \times^7 S^1$ such that $r^2 = (tuv)/(TUV)$.  We can send a   $\vek{c}$
  coordinate  to such a coordinates via
\[ \vek{b}(\vek{c}) =  
\left(\sr{A}\sr{D},\sr{B}\sr{E},\sr{C}\sr{F},\sr{D}/\sr{A},\sr{E}/\sr{B},\sr{F}/\sr{C},- \sr{A}\sr{B}\sr{C}
\right),
\] 
and we can send  a  $\vek{b}$ coordinate to a $\vek{s}$  coordinate 
via 
\[  \vek{s}(\vek{b}) =(t/T,u/U,v/V,tT,uU,vV,UVr,TUr TVr, r). \]
As such, we find that the  
 the tetrahedral circulants cover the balanced coordinates via 4 fold cover  (as described in the comment at the end of section \ref{com}) and the balance coordinates cover the $\vek{s}$  coordinates via a 2 fold cover, given by  negating $t,u,v,T,U$ and $V$ and leaving $r$ alone.  We will call this transformation $D$, 
 and note that 
 $\vek{b}(\hat{D} \cdot \vek{c})=  D \cdot \vek{b}(\vek{c})$, when $\hat{D} \in \negg$.   From this identity,  the 
balanced  coordinates parametrize our \op.

Using the balanced coordinates, we provide an explicit algebraic description of the group $\wacky$ from  section \ref{vsym}.
To do so,  we first lift the action on the $\vek{s}$ coordinate described by  equations (\ref{g1})-(\ref{g12})    to a transformation of  the balanced coordinates. We will index our group elements  by how they act on $(t,u,v,T,U,V;r)$.
Using the same argument as found in lemma \ref{inter2}, we can then determine whether our lift or $D$ times our lift represents the needed transformation.   We find  the shaded elements are given by 
\[   \left\{ [\bar{T},u,v,\bar{t},U,V;  r],
[T,u,v,t,U,V;  T \bar{t} r],[t,U,v,T,u,V;  U \bar{u} r] \right. \]
\[\left. ,[t,u,V,T,U,v;  V \bar{v} r],
[t,\bar{U},v,T,\bar{u},V;  r], [t,u,\bar{V},T,U,\bar{v};  r]   
\right\}. \]
while the unshaded elements are given by 
\[ \left\{  [U,u,v,T,t,V;  U \bar{t} r],[\bar{U},u,v,T,\bar{t},V ; r] , [t,V,v,T,U,u;   V \bar{u} r], 
\right. \]
\[ \left. [t,\bar{V},v,T,U,\bar{u}; \bar{u}\bar{V} r] , [t,u,T,v,U,V;   T \bar{v} r] ,[t,u,\bar{T},\bar{v},U,V; \bar{v}\bar{T} r]\right\}. \]
Note the action on $r$ is determined by the action on the first 6 coordinates,
and that $r$ does not affect the action upon the first six coordinates. 
Hence to identify the group we only need to understand $\wacky$'s action on the 
$(T,u,v,t,U,V)$ coordinates. 
 Each of these generators   is a  permutations of the  $\vek{b}$ coordinates together with  something in the group generated by conjugating an even number of the $\vek{b}$ coordinates.      In fact,  these elements are easily checked to generate this group.  This is a well known reflection group   usually  denoted as $D_6$.

{\bf Note:} In the $\vek{b}$ coordinates
\[ \alpha(\vek{b})  =2 ( t+ u + v +  tuv - r (T+U+V+TUV)),\]
\[ \beta(\vek{b})  = (T+U+V+1/T+1/U+1/V)-(t+u+v+1/t+1/u+1/v).\]
 \[ \vec{\gamma}(\vek{b}) =[r,1,rTU,tu,rTV,tv,rUV,uv] .\] 
 Due to the simple nature of $\wacky$ in the balanced coordinates, it is easy to explore $\wacky$'s action on these quantities.  For example, 
by  plugging in the above generators, we find that  $\delta$ is $\wacky$ invariant.  


\end{subsubsection}

\begin{subsubsection}{Scissors Cosets  \label{sispre}}

With the balanced coordinates we can easily   describe our  30 scissors classes described in the proof of  corollary \ref{cor}. From equation \ref{mirror}, the scissor classes described in corollary \ref{cor} occur in pairs consisting of a tetrahedron and its mirror image.  
Hence, we only need to describe fifteen  scissors classes no pair of which are mirror images of each other. In other words, we need to describe  the cosets of the quotient of $\wacky$ by the group generated by the tetrahedral symmetries together with the shaded subgroup.  To do so, we let a lower and capital case coordinate pair in $(t,u,v,T,U,V)$, like $t$ and $T$,  be called a {\it pair}.  The shade subgroup in balanced coordinates  is generated by the transformations that swap the elements of a pair, for example $TuvtUV$, and the elements that  conjugate a pair, for example $\bar{t}uv\bar{T}UV$.  The group of tetrahedral symmetries is generated by permuting the pairs, for example $utvUTV$,  together with the transformations that conjugate a pair of the capitol letters, for example
 $tuv\bar{T}\bar{U}V$.  
Together, the shaded subgroup and tetrahedral symmetries  form a group  generated by  all pair swaps,  all pair  permutations, together with all even conjugations.   Hence 
the needed cosets are indexed by the elements of the following sets:
\[ SR = \{  tuvTUV ,tuvTVU,
 tuvUTV, tuvUVT,
 tuvVTU,tuvVUT \}, \] 
together with  
\[ SN =\left\{  tuTvVU, tuTvUV  ,tuUvTV,
tvTuVU,tvTuUV, \ tvUuTV, \right. \]
\[ \left.  uvTtVU,uvTtUV,uvUtTV,
\right\}, \]
where the group elements have been indexed by how they  act on $(t,u,v,T,U,V)$.
We distinguished between these two subsets since the $SR$ corresponds to the 6 nontrivial Regge scissors classes, which, in the $\vek{b}$ coordinates, is the group generated by independently permuting the lower or upper case coordinates together with the even conjugations of the uppercase coordinates.   $SN$ 
corresponds to the remaining 9 nontrivial  scissors classes. 

\end{subsubsection}

\end{subsection}

\begin{subsection}{The Generic Set \label{dis}}
 
 We will now describe an open dense set of the \op where all the construction needed to prove theorems  \ref{fundsis} and \ref{res2}  are guaranteed to apply. 
 Namely, we will restrict the $\gen \cdot \vek{b}$ to  where  the 8 ideal tetrahedra from figure \ref{fireoct}  and the 6 
 ideal tetrahedra arising when inverting $(\gen \cdot \vek{b})_{oct}$ to a supertetrahedron,   
 as in figure \ref{heart}, are non-degenerate for all $\gen \in G$.  When all such tetrahedra are nondegenerate, all our construction make sense with out modification. 
 
 \begin{lemma}\label{sil1} 
Suppose all  images of $r$ under $\wacky$ are nondegenerate,
then the ideal tetrahedra forming $(\gen \cdot \vek{b})_{oct}$,   as in figure \ref{fireoct}, are 
 nondegenerate for every $g  \in G$.
 \end{lemma}
\proof 
We need to check that all the tetrahedral clinants of $(\vek{b})_{oct}$  are nondegenerate under this assumption. Note that the $\wacky$ images of $r$ under $\wacky$ are the `positive' square roots of all the even conjugations applied to $tuvTUV$.  
  The waist has clinants which are all in this form.  To see that the remaining tetrahedral clinants are nondegenerate,  note that the $z$ from figure \ref{fireoct} satisfies equation \ref{hol2}, and, 
  hence,
    \begin{eqnarray}\label{hol3}
 \Pi_{i=1}^4(1-  \overline{(\vek{s}(\vek{b}),z)^{\mat{}}_{\sr{M}\sr{Y}}(2i)} )  & = &   \Pi_{i=1}^4(1-  (\vek{s}(\vek{b}),z)^{\mat{}}_{\sr{M}\sr{Y}}(2i-1)  ) .
   \end{eqnarray}
  Notice the $(\vek{s}(\vek{b}),z)^{\mat{}}_{\sr{M}\sr{Y}}(j)$ are our need tetrahedral clinants, up to conjugation. 
Let us assume one of these clinants is degenerate and produce a contradiction.  Under this assumption to one side of equation (\ref{hol3}) is equal to zero.  If one side is zero then the other has to be zero, hence, have a  degenerate term as well.  So, for some $i$ and $j$  
 \begin{eqnarray}\label{trap}
 (\overline{(\vek{s}(\vek{b}),z)^{\mat{}}_{\sr{M}\sr{Y}}(2i)} )(\vek{s}(\vek{b}),z)^{\mat{}}_{\sr{M}\sr{Y}}(2j-1) =1.
 \end{eqnarray}
 Upon multiply this expression out, we find that the right hand side of equation (\ref{trap}) is 
 an  image of $r$ under $\wacky$,  our needed contradiction.
 
 \qed

This lemma assures us that when  acting on the \op that we never run into a degeneracy in the compliment of the set determined by the  16  monomial constraints derived from $r=1$ under $\wacky$'s action.   Notice these clinants correspond to supertetrahedral  dihedral clinants under the edge flips alone.  Hence, none are found in the interior of the collection of the \op  corresponding to the finite tetrahedra and all there images under $\wacky$.   This is a useful observation, for it tells us that if we look at the image of the space  finite tetrahedra under 
$\gen \in G$, then   whether a triple of planes intersects at a finite or hyperideal point is the same for the image of every finite tetrahedron.  Clearly the  group generated by the shaded subgroup and the tetrahedral symmetries satisfy that such intersections occur at finite vertices.  Hence, we can understand every element of $\wacky$ by examining the 15 cosets in section \ref{sispre}.  We find that the 6 Regge cosets in $SR$ preserve the fact that all the planes intersect at finite points, while the remaining 9 cosets have all 4 triples of planes intersecting at hyperideal points.

 \begin{lemma}\label{sil2}
Suppose  no clinant  derived from  $r$ or $tT$ via the action of $\wacky$ is degenerate, 
then the ideal tetrahedra utilized when forming the supertetrahedron from $(\gen \cdot \vek{b})_{oct}$,   as in figure \ref{heart}, are 
 nondegenerate for every $g \in G$.
 \end{lemma}
\proof
The clinants of these tetrahedron are either an  image of $r$ under $\wacky$ or supertetrahedral   clinant along a tetrahedral edge.  By the definition of the balance coordinates the supertetrahedral clinants along a tetrahedral edge
 are in the form  of  an element of $\wacky$ applied to  $tT$, as needed.   

\qed

From lemmas \ref{sil1} and \ref{sil2},  we will never run into a degenerate tetrahedron when working 
in  the set where $(t,u,v,T,U,V)$ is  in the compliment of the set determined by the  46  monomial constraints derived from  $tT=1$ and $r=1$
under $\wacky$'s action.   
   
  {\bf Comment:}   In relatively straight-forward ways, the constructions presented in section \ref{simp} and \ref{clis}  can be extended to much of 
  the set where an image of $r$ or $tT$ under the action of $\wacky$ is degenerate.  However, attempts at finding a  unified approach have not been successful.     In order to make good geometric sense out of these constructions in general, we must take the correct compactification of the set of $\vek{b}$ where our constraints are  nondegenerate.  This will force us   to  go beyond the $C$-regions, see the note at the end of section \ref{simp}.  The above constraints eliminate many interesting cases, like a tetrahedron with ideal vertices
 and any  tetrahedron exhibiting a symmetry 
where, for example,  $\sr{A} \sr{B} = \sr{D} \sr{E}$.


 \end{subsection}
 
 \end{section}
 
 \begin{section}{Derivation of Equation (\ref{nicev}) \label{dnicev}}

Recall from section \ref{alt}, that to derive equation (\ref{nicev}), we first derive equation (\ref{hnice}), which we do in section \ref{dhnice}.  Then we check that $\scr{H}(w)$ is odd, which is verified in section \ref{FSymd}.

  \begin{subsection}{Derivation of Equation (\ref{hnice}) \label{dhnice}}   
   
\begin{figure}
\center{\includegraphics{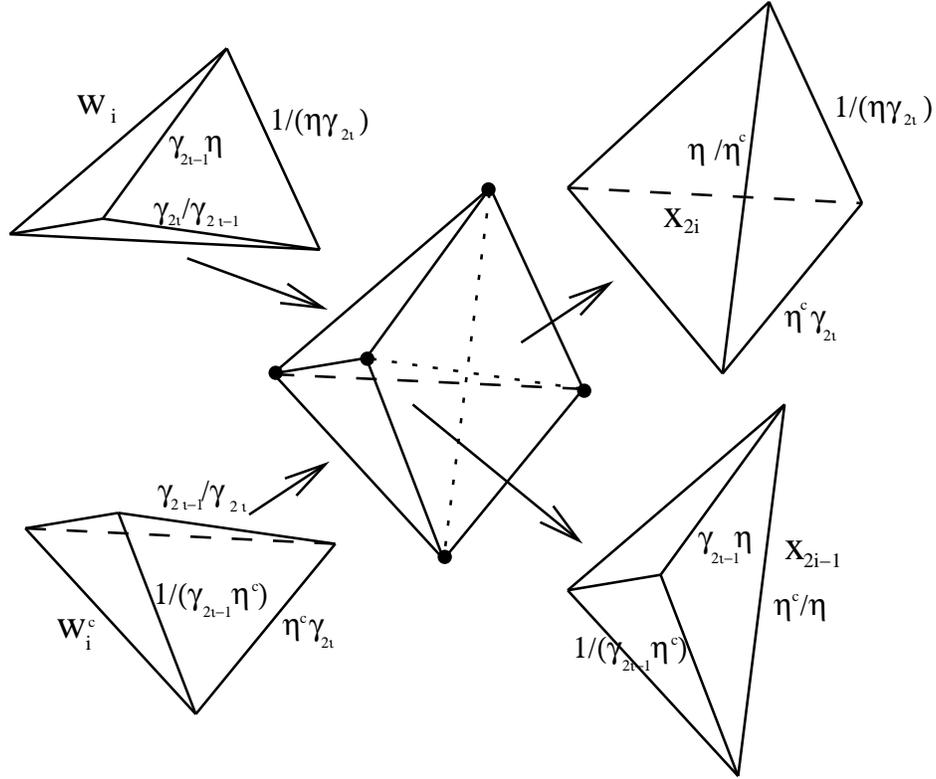}}
\caption{
\label{deggy}
On the left hand side of this figure we have a pair of the 
ideal   tetrahedra arising in the description of 
$ (\vek{b})_{op} $. 
 To be precise, we have pictured  
$w_i(\vek{b} )$, labeled  as $w_i$,  and  $w_i(\bar{\vek{b}})$,  labeled as  $w_i^c$.
We have labeled the tetrahedral clients using the $\vec{\gamma}$ from section \ref{alt}, evaluated at  $\osa^{-1} \cdot \vek{b}$, 
or rather $\vec{\gamma} (\osa^{-1} \cdot \vek{b}) =   \left(
1,
\frac{1}{ r_2},
\frac{1}{b},
\frac{c}{r_4},
\frac{c}{ r_3 r_4},
\frac{1}{ r_4},
\frac{   r_1}{ a r_2}
\frac{f}{ r_2}
 \right).
$
We glue our   $w_i$ and $w_i^c$ together in order to form the $C$-region in the middle of this picture.
We  then resplit this $C$-region to construct the labeled
$x_{j}$ complex coordinates.  We let $z_{2i-1}= 1/x_{2i-1}$ and  $z_{2i}= x_{2i}$.    Let us denote  the list of 8 ideal tetrahedra indexed by these $z_i$ as  $\vek{z} \in {\bf C}^8 $. 
We chose to let $z_{2i-1}$  to be one over the $x_{2i}$ coordinate of the tetrahedron  in the lower right corner of our figure so that all the $z_j$ share the same  clinant, $\eta/\eta^c$, where $\eta=\mat{\rho}(\vek{b})$ and $\eta^c=\overline{\mat{\rho}(\bar{\vek{b}})}$. }
\end{figure}

The first step in proving equation  (\ref{hnice})  is the geometric rearrangement of  \op as described in figure \ref{deggy}.       In this section,  we will assume $\alpha$, $\beta$,  and $\gamma$ are evaluated at $\osa^{-1} \cdot \vek{b}$ which is equivalent to the hats in section \ref{clis}.
From figure \ref{deggy}, we have that 
\begin{eqnarray}\label{zform}
2  [(\vek{b} )_{op}] & =  & \sum_{i=1}^8 (-1)^{i}  [z_i (\vek{b})] .
 \end{eqnarray}
and, utilizing equation (\ref{fund}), we have  
 \[ z_{j}(\vek{b}) = \frac{1- \gamma_{j}  \eta}{1-  \gamma_{j}  \eta^c}. \]
The $\vek{z}$ form coordinates on  the space of  \op
that  reside in the following set:
\[  \scr{DT} =  \left\{ \vek{z} \in {\bf C}^8  \left|   \Pi_i^4 \frac{z_{2i-1}}{z_{2i}} =1, \Pi_i^4 \frac{1-z_{2i-1}}{1-z_{2i}} =1,Im\left( \frac{z_i}{z_j} \right) = 0 \right. \right\}.\]
Conversely, by inverting the construction in figure \ref{deggy}, the elements of $\scr{DT}$ can be used to form \op.
Hence, the $\vek{z}$ in $\scr{DT}$  parameterize  the \op.

Given a complex number $z$ we let $z/|z|$ be called its associated {\it circulant} and  $z/\bar{z}$ be called  its associated {\it clinant}.
The condition $ \Im\left( \frac{z_i}{z_j} \right) = 0$ is equivalent to the fact that all our $z_i$ have the same clinant associated to them, 
$\eta/\eta^c$  as seen in figure \ref{deggy}. We will call this  the {\it magic clinant}.  The magic clinant is preserved by an order-2304 subgroup $H < G$ described in section \ref{revis}. 
Representatives of the ten right cosets of $H$ in $\wacky$ index ten truly distinct volume formulas.
These cosets are explicitly described in section \ref{revis}.

Since every angle  of an octahedron is on the waist with respect  to some simplicial decomposition,  the octahedron, unlike our other $C$-regions,  has well defined circulants. 
From figure \ref{deggy}, we know that  the $z_i$ share the same clinant. However, the $z_i$ will not, in general,  share the same circulant.
For  convex supertetrahedron,  by explicitly performing the construction in figure \ref{heart},  we find that the even and odd $z_i$  have circulants that differ by a minus sign.
In fact, both the \op corresponding to a convex supertetrahedron are convex, have the same orientation  and one has dihedral angles that are $\pi$ minus the other's dihedral angles.  
As another example,  when we cut down a standard finite vertexed supertetrahedron we find that the resulting $z_i$ have circulants that all agree.
In fact, by explicitly performing the construction in figure \ref{heart}, one of the \op corresponding to a  standard finite vertexed supertetrahedron 
is the opposite orientation of the other and has negated dihedral angles. 
Hence, the finite tetrahedra correspond to a  subset of 
\[ \scr{FT} = \left\{ \vek{z} \in {\bf C}^8  \left|   \Pi_i^4 \frac{z_{2i-1}}{z_{2i}} =1, \Pi_i^4 \frac{1-z_{2i-1}}{1-z_{2i}} =1, \frac{z_1}{z_j} > 0 \right. \right\}.\]

We now derive formula (\ref{hnice}).
Notice that from equation (\ref{cross}), that 
\[ \scr{B}(z) = -\scr{B}(1-z) = -\Im(\scr{L}_2(1-z)  + I\arg(z) \log|1-z|). \]
Hence, from  
equations (\ref{zform}) and (\ref{iv1}), we have that 
\[  \vol(2[(\vek{b})_{tet}] ) =  -\sum_{i=1}^8 (-1)^{i} \scr{B}(1-z_i). \]
If we restrict our selves to  $\scr{FT}$, then all the $z_i$ circulants are the same, and  the fact that 
\[ \Pi_i^4 \left| \frac{1-z_{2i-1}}{1-z_{2i}} \right| = 1\]
assures us that we can ignore  the  $\arg(z) \log|1-z|$  part of $\scr{B}$.
Hence we have 
\[  \vol(2[(\vek{b})_{tet}] ) =- \Im \left( \sum_{i=1}^8 (-1)^{i} \scr{L}_2(1-z_i) \right). \]
For $\vek{z}$ in $\scr{DT}- \scr{FT}$ the  circulants don't agree and we will have to tack on $\pm I \pi \log(z)$ terms.  For simplicity we will derive our formulas in the $\scr{FT}$ case.   
In this case, the fact  that 
\[ \Pi_i^4 \left| \frac{z_{2i-1}}{z_{2i}} \right| = 1.\]
assures us that 
\[ \Im \left( \sum_{i=1}^8 (-1)^{i} \log(z_i)^2 \right)=0. \]
(we use the $\log(z)$ with $\log(1)=0$ and a branch cut along $(-\infty,0]$).
Hence, 
 \[ \scr{L}(z) =\scr{L}_2(1-z)+  \frac{1}{4}\log(z)^2 ,  \] 
 is an  analytic function (with a branch cut along $(-\infty,0]$  and  $\scr{L}(1)=0$) such  that 
\[  \vol(2[(\vek{b})_{tet}] )  =  \Im( \sum_{i=1}^8 (-1)^{i} \scr{L}(z_i)). \]
Now we let 
\[ \scr{H}(w) = \scr{L}\left(  \frac{1-w}{1+w} \right), \]
with branch cuts  $[-\infty,1] \bigcup [1,\infty]$ and $H(0)=0$. 
Note
$ \frac{1-w}{1+w}$ is its own inverse and 
\begin{eqnarray*}\label{FSymmy}
\frac{1-z_i}{1+z_i}& =  &  \frac{{\gamma}_i(\eta-\eta^c)}{2-{\gamma_i} (\eta+\eta^c)}  \\
& = & \frac{-\gamma_i(2 i \sqrt{\delta}/{\alpha})}{2+2 {\gamma}_i( {\beta}/ {\alpha})}    \\
 & =  &  \frac{-i  {\gamma}_i  \sqrt{\delta}  }{ {\alpha}+  {\gamma}_i( {\beta})}   \\
  & =  &  \sqrt{\delta}  {c}_i  
 \end{eqnarray*}
 Where 
 \[ {c}_i(\osa^{-1} \cdot \vek{b}) = \frac{-i  {\gamma}_i}{{{\alpha}}+{\gamma}_i {\beta}}. \] 
 Hence 
  \begin{eqnarray*}\label{FSym1}
\scr{L}\left(  z_i \right) & =  & \scr{H}\left(  \frac{1-z_i}{1+z_i} \right) \\
& = &  \scr{H}\left( {c}_i \sqrt{\delta} \right) 
 \end{eqnarray*}

 We nearly have 
equation (\ref{nicev}).  Namely, everything has been evaluated at 
$\osa^{-1} \cdot \vek{b}$, hence we must apply $\osa$.
A priori, there is a possibility that we will need to tack on tack on $\pm I \pi \log(z)$ terms to our $\scr{L}$.
To assures ourselves that we do not need to,  we must check that the $z_i/z_j >0$ constraint is preserved under the $\osa$ transformation.
Such a constraint can only changes signs if an ideal  tetrahedra in $(\vek{b})_{op}$  degenerates.
By lemma \ref{sil1}, this cannot occur in the connected  set of finite tetrahedra, hence cannot occur in any image of this set under $\wacky$.  
So we only need to verify this for a single finite  tetrahedron, which is easily accomplished.

{\bf Comment 1:} In this section, we see  that  we only need the  generalized tetrahedra in $\scr{FT} \subset \scr{DT}$ in order to derive equation (\ref{MYF}).  The finite tetrahedra are a proper subset of $\scr{FT}$ and generalized tetrahedra are still necessary.  For example, the  transformation $\osa$ will send any finite hyperbolic  tetrahedron to a tetrahedron with tetrahedral planes that intersect at hyperideal points.     

{\bf Comment 2:} The $\scr{DT}$ parameterization of our \op  makes transparent many of our $\wacky$ symmetries.  Namely,
$\scr{DT}$ is clearly preserved by independently permuting  the even and odd indexed $z_i$, by  the transformation that sends the $z$ forming $\vek{z}$ to  
 $\frac{1}{\bar{z}}$, and  by the transformation that  conjugates and swaps the $z_{2i-1}$  and $z_{2i}$ terms.  These transformations clearly preserve scissors class and, in fact, generate the    $H$ subgroup of $\wacky$ discussed  
 in section \ref{revis}.

\end{subsection}

 \begin{subsection}{$\scr{H}(w)$ is odd  \label{FSymd}}

 Notice that $H(w)$ is odd if 
 $\scr{L}(z)$  from section \ref{dhnice}  satisfies 
  \begin{eqnarray} \label{FSyml}
\scr{L} \left(\frac{1}{z}\right) & =  & -\scr{L} ( z ).
 \end{eqnarray}
In order to prove equation (\ref{FSyml}),
first  we  define 
\begin{eqnarray}\label{ob}
\scr{K}(z) & =  &  \frac{\scr{L}_{2}(1-z) - \scr{L}_2(1-1/z)}{2}  
\end{eqnarray}
which clearly satisfies
  \begin{eqnarray} \label{FSymk}
\scr{K} \left(\frac{1}{z}\right) & =  & -\scr{K} ( z ).
 \end{eqnarray}
 Hence,   equation (\ref{FSyml}) will follow from 
 the following lemma.
 \begin{lemma}
\[    \scr{L}(z)= \scr{K}(z) \]
\end{lemma}
\proof  
This lemma will follows if we can demonstrate that  
\begin{eqnarray}\label{FSym2}
  \scr{L}_2(1-1/z) &  = & - \scr{L}_{2}(1-z) - \frac{1}{2} \log(z)^2 .
  \end{eqnarray}
To prove equation (\ref{FSym2}), 
we rewrite  equation (\ref{FSym2}) as 
\begin{eqnarray}\label{FSym3}
  \scr{L}_2(1-1/z)  + \scr{L}_{2}(1-z)&  = & - \frac{1}{2} \log(z)^2 
  \end{eqnarray}
and note that the right and left hand sides of  equation (\ref{FSym3}) are equal to 0 at $z=1$.  Hence we need only show that 
\[  d\left(\scr{L}_2(1-z) + \scr{L}_2(1-1/z) \right) =  \frac{-1}{2} d \log(z)^2 \] 
We compute the left hand side
\begin{eqnarray*}
  d(\scr{L}_2(1-z) + \scr{L}_2(1-\frac{1}{z})) & = & -\log(z)d(\log(1-z)) -\log\left(1/z \right)
  d(\log(1-1/z) ) \\
   & = &  -\log(z)d(\log(1-z)) + \log(z)d\left(\log\left((z-1)/z \right)\right)   \\
    & = & - \log(z)d(\log(z))       
    \end{eqnarray*}
and find it is   $\frac{-1}{2} d \log(z)^2$ as needed. 

\qed

{\bf Comment:} Equation (\ref{FSyml}) is   very suggestive with regards to the Chern Simons invariant, but this turns out to be misleading since the scissors congruence taking $\vek{z}$ to $\vek{\frac{1}{\bar{z}}}$,  requires using the fact 
  that  $[z]=-[\bar{z}] $, which destroys the type of orientation-sensitivity needed to capture the 
  Chern Simons invariant (see Neumann \cite{Ne}).

\end{subsection}

\begin{subsection}{10 Interesting Formulas  \label{revis}}

Once one has a volume formula in hand one can re-write this formula in many ways utilizing $\wacky$, as described in  section \ref{algsym}. 
For example we have that 
\[ \vol([(\vek{b})_{tet}]) =
 \sum_{i=1}^8 (-1)^{i} 
\scr{L}( z_i(g \cdot (\vek{b})_{tet}) ) \]
for every $\gen \in G$. 
What really distinguishes these different formulas is the  magic clinant, as introduced in section \ref{dhnice}.   The magic clinant can be easily explored in the $\vek{b}$ coordinates by noting 
\[  \eta(\vek{b})/\eta(\vek{b})^c =  \left( \frac{-\beta+ i \sqrt{\delta}}{\alpha} \right) 
 \left( \frac{-\beta + i \sqrt{\delta}}{\bar{\alpha}} \right) =  \frac{(\beta^2- \delta)-i (2 \beta \sqrt{\delta})}{|\alpha|^2}   \]
 or rather
 \begin{eqnarray}\label{magcli}
 m(\vek{b}) & = & \frac{(\beta^2- \delta)-i (2 \beta \sqrt{\delta})}{|\alpha|^2}. 
 \end{eqnarray}
In particular, we can  explore the subgroup of $\wacky$  that preserves the magic clinant.
Let $H_0$ be the subgroup  generated by  independently permuting the lower and upper case coordinates and  performing an even number of conjugations.
Since, as in the note at the end of section \ref{algsym}, all $\gen \in G$ preserve $\delta$, and $H_0$ clearly  preserve $\beta$, we have that $H_0$ must  must preserve $|\alpha|^2$. Hence, from equation (\ref{magcli}), $H_0$ preserves the magic clinant. 
Similarly, the transformation that   swaps all  pairs simultaneously  (as defined in section \ref{sispre}) will    negate $\beta$.  Hence, this transformation conjugates the magic clinant. 
 Let $H$ is be the order-2304 subgroup generated by this swap transformation and $H_0$. The elements of $\wacky$ not in $H$ will not preserve the magic clinant.  From this 
observation, we arrive at 10 truly distinct volume formulas index by $G/H$, with coset representatives given by 
\[
\{ 
tuvTUV ,TuvtUV, t UvTuV,tuVTUv, tvVTuU \] \[ ,tTVuUv,tTuUvV 
, tTUuvV ,tuUTvV,tTvuUV \} .
\] 



\end{subsection}

 \end{section}
 
\begin{section}{Questions}

 
{\bf Question 1:}  Numerically, equation (\ref{MYF}) holds in the spherical case and, by the analytic
continuation principle, this is not very surprising.   However, one finds that some of the 
constructions presented here become difficult to implement in the  spherical case.  Can these constructions be
made to make sense in the spherical world?  In particular, can one prove theorem \ref{fundsis} in the
spherical case? \vspace{.1in}  \linebreak 
 {\bf Question 2:} Notice the scissor group produced here has been explicitly described 
when it is acting between pairs of tetrahedra.  
In order reduce to the tetrahedra, we are forced to use equation 
 \ref{half}.   
 Dupont and Sah's    division algorithm, utilized  to prove equation 
 \ref{half},  is rather complicated to implement geometrically  (see
\cite{Sa1}).  Is there a simple way to accomplish this
division in this case?

\end{section}

\vspace{.2in} 

 {\bf Acknowledgments:} 
 The authors like  to thank Yana  Mohanty, Dylan Thurston, and Walter Neumann for the useful discussions we had with them concerning this work.

 \bibliographystyle{plain}
\bibliography{VSbib}

\end{document}